\documentclass{article}%
\usepackage{amsfonts}
\usepackage{amsmath}
\usepackage{amssymb}
\usepackage{graphicx}%
\newtheorem{theorem}{Theorem}

\newtheorem{definition}[theorem]{Definition}

\newtheorem{lemma}[theorem]{Lemma}

\newtheorem{remark}[theorem]{Remark}

\newenvironment{proof}[1][Proof]{\noindent\textbf{#1.} }{\ \rule{0.5em}{0.5em}}
\begin{document}

\title{Necessary and sufficient optimality conditions for relaxed and strict control
problems of backward systems}
\author{\textbf{Seid BAHLALI}\\Laboratory of Applied Mathematics, University Med Khider,\\Po. Box 145, Biskra 07000, Algeria.\\sbahlali@yahoo.fr \ }
\maketitle

\begin{abstract}
We consider a stochastic control problem where the set of strict (classical)
controls is not necessarily convex, and the system is governed by a nonlinear
backward stochastic differential equation. By introducing a new approach, we
establish necessary as well as sufficient conditions of optimality for two
models. The first concerns the relaxed controls, who are measure-valued
processes. The second is a particular case of the first and relates to strict
control problems.

\ 


\ 


\end{abstract}

\section{Introduction}

We study a stochastic control problem where the system is governed by a
nonlinear backward stochastic differential equation (BSDE\ for short) of the
type
\[
\left\{
\begin{array}
[c]{l}%
dy_{t}=b\left(  t,y_{t}^{v},z_{t}^{v},v_{t}\right)  dt+z_{t}^{v}dW_{t},\\
y_{T}=\xi
\end{array}
\right.
\]
where $b$ is given maps, $W=\left(  W_{t}\right)  _{t\geq0}$ is a standard
Brownian motion, defined on a filtered probability space $\left(
\Omega,\mathcal{F},\left(  \mathcal{F}_{t}\right)  _{t\geq0}%
,\mathbb{\mathcal{P}}\right)  ,$ satisfying the usual conditions.

The control variable $v=\left(  v_{t}\right)  $, called strict (classical)
control, is an $\mathcal{F}_{t}$ adapted process with values in some set $U$
of $\mathbb{R}^{m}$. We denote by $\mathcal{U}$ the class of all strict controls.

The criteria to be minimized, over the set $\mathcal{U}$, has the form%
\[
J\left(  v\right)  =\mathbb{E}\left[  g\left(  y_{0}^{v}\right)  +%
{\displaystyle\int\nolimits_{0}^{T}}
h\left(  t,y_{t}^{v},z_{t}^{v},v_{t}\right)  dt\right]  ,
\]
where $g$ and$\ h$ are given functions and $\left(  y_{t}^{v},z_{t}%
^{v}\right)  $ is the trajectory of the system controlled by $v.$

A control $u\in\mathcal{U}$ is called optimal if it satisfies%
\[
J\left(  u\right)  =\underset{v\in\mathcal{U}}{\inf}J\left(  v\right)  .
\]

Stochastic control problems for backward and forward-backward systems have
been studied by many authors including Peng $\left[  27\right]  $, Xu $\left[
31\right]  $, El-Karoui et al $\left[  12,13\right]  $, Wu $\left[  30\right]
$, Dokuchaev and Zhou $\left[  9\right]  $, Peng and Wu $\left[  28\right]  $,
Bahlali and Labed $\left[  2\right]  $, Bahlali $\left[  5,6\right]  $, Shi
and Wu $\left[  29\right]  $, Ji and Zhou $\left[  19\right]  $. The dynamic
programming approach was studied by Fuhrman and Tessitore $\left[  16\right]
.$

Since the strict control domain being nonconvex, then if we use the classical
method of spike variation on strict controls, the major difficulty in doing
this is that the generator $b$ and the running cost coefficient $h$ depend on
two variables $y_{t}$\ and $z_{t}$. Then, we can't derive directly the
variational inequality, because $z_{t}$\ is hard to handle, there is no
convenient pointwise (in $t$) estimation for it, as opposed to the first
variable $y_{t}$. To overcome this difficulty, we introduce a new approach
which consist to use a bigger new class $\mathcal{R}$ of processes by
replacing the $U$-valued process $\left(  v_{t}\right)  $ by a $\mathbb{P}%
\left(  U\right)  $-valued process $\left(  q_{t}\right)  $, where
$\mathbb{P}\left(  U\right)  $ is the space of probability measures on $U$
equipped with the topology of stable convergence. This new class of processes
is called relaxed controls and have a richer structure of compacity and
convexity. This property of convexity of relaxed controls, enables us to treat
the problem with the way of convex perturbation on relaxed controls.

In the relaxed model, the system is governed by the BSDE%
\[
\left\{
\begin{array}
[c]{l}%
dy_{t}^{q}=\int_{U}b\left(  t,y_{t}^{q},z_{t}^{q},a\right)  q_{t}\left(
da\right)  dt+z_{t}^{q}dW_{t},\\
y_{T}^{q}=\xi.
\end{array}
\right.
\]

The functional cost to be minimized, over the class $\mathcal{R}$ of relaxed
controls, is defined by%
\[
J\left(  q\right)  =\mathbb{E}\left[  g\left(  y_{0}^{q}\right)  +%
{\displaystyle\int\nolimits_{0}^{T}}
\int_{U}h\left(  t,y_{t}^{q},z_{t}^{q},a\right)  q_{t}\left(  da\right)
dt\right]  .
\]

A relaxed control $\mu$ is called optimal if it solves
\[
J\left(  \mu\right)  =\inf\limits_{q\in\mathcal{R}}J\left(  q\right)  .
\]

The relaxed control problem is a generalization of the problem of strict
controls. Indeed, if $q_{t}\left(  da\right)  =\delta_{v_{t}}\left(
da\right)  $ is a Dirac measure concentrated at a single point $v_{t}\in U$,
then we get a strict control problem as a particular case of the relaxed one.

Our aim in this paper, is to establish necessary as well as sufficient
conditions of optimality in the form of global stochastic maximum principle,
for both relaxed and strict controls. To achieve this goal, we derive these
results as follows.

Firstly, we give the optimality conditions for relaxed controls. The idea is
to use the fact that the set of relaxed controls is convex. Then, we establish
necessary optimality conditions by using the classical way of the convex
perturbation method. More precisely, if we denote by $\mu$ an optimal relaxed
control and $q$ is an arbitrary element of $\mathcal{R}$, then with a
sufficiently small $\theta>0$ and for each $t\in\left[  0,T\right]  $, we can
define a perturbed control as follows%
\[
\mu_{t}^{\theta}=\mu_{t}+\theta\left(  q_{t}-\mu_{t}\right)  .
\]

We derive the variational equation from the state equation, and the
variational inequality from the inequality%
\[
0\leq J\left(  \mu^{\theta}\right)  -J\left(  \mu\right)  .
\]

By using the fact that the coefficients $b$ and $h$ are linear with respect to
the relaxed control variable, necessary optimality conditions are obtained
directly in the global form.

To achieve this part of the paper, we prove under minimal additional
hypothesis, that these necessary optimality conditions for relaxed controls
are also sufficient.

The second main result in the paper characterizes the optimality for strict
control processes. It is directly derived from the above result by restricting
from relaxed to strict controls. The idea is to replace the relaxed controls
by a Dirac measures charging a strict controls. Thus, we reduce the set
$\mathcal{R}$ of relaxed controls and we minimize the cost $J$ over the subset
$\delta\left(  \mathcal{U}\right)  =\left\{  q\in\mathcal{R}\text{ \ /
}\ q=\delta_{v}\ \ ;\ \ v\in\mathcal{U}\right\}  $. Necessary optimality
conditions for strict controls are then obtained directly from those of
relaxed one. Finally, we prove that these necessary conditions becomes
sufficient, without imposing neither the convexity of $U$ nor that of the
Hamiltonian $H$ in $v$.

The paper is organized as follows. In Section 2, we formulate the strict and
relaxed control problems and give the various assumptions used throughout the
paper. Section 3 is devoted to study the relaxed control problems and we
establish necessary as well as sufficient conditions of optimality for relaxed
controls. In the last section, we derive directly from the results of Section
3, the optimality conditions for strict controls.

\ 

Along this paper, we denote by $C$ some positive constant and we need the
following matrix notations. We denote by $\mathcal{M}_{n\times d}\left(
\mathbb{R}\right)  $ the space of $n\times d$ real matrices and by
$\mathcal{M}_{n\times n}^{d}\left(  \mathbb{R}\right)  $ the linear space of
vectors $M=\left(  M_{1},...,M_{d}\right)  $ where $M_{i}\in\mathcal{M}%
_{n\times n}\left(  \mathbb{R}\right)  $.

For any $M,N\in\mathcal{M}_{n\times n}^{d}\left(  \mathbb{R}\right)  $,
$L,S\in\mathcal{M}_{n\times d}\left(  \mathbb{R}\right)  $, $Q\in
\mathcal{M}_{n\times n}\left(  \mathbb{R}\right)  $, $\alpha,\beta
\in\mathbb{R}^{n}$ and $\gamma\in\mathbb{R}^{d},$ we use the following notations

$\alpha\beta=%
{\displaystyle\sum\limits_{i=1}^{n}}
\alpha_{i}\beta_{i}\in\mathbb{R}$ is the product scalar in $\mathbb{R}^{n}$;

$LS=%
{\displaystyle\sum\limits_{i=1}^{d}}
L_{i}S_{i}\in\mathbb{R}$, where $L_{i}$ and\ $S_{i}$ are the $i^{th}$ columns
of $L$ and $S;$

$ML=%
{\displaystyle\sum\limits_{i=1}^{d}}
M_{i}L_{i}\in\mathbb{R}^{n}$;

$M\alpha\gamma=\sum\limits_{i=1}^{d}\left(  M_{i}\alpha\right)  \gamma_{i}%
\in\mathbb{R}^{n}$;

$MN=%
{\displaystyle\sum\limits_{i=1}^{d}}
M_{i}N_{i}\in\mathcal{M}_{n\times n}\left(  \mathbb{R}\right)  $;

$MQN=%
{\displaystyle\sum\limits_{i=1}^{d}}
M_{i}QN_{i}\in\mathcal{M}_{n\times n}\left(  \mathbb{R}\right)  $;

$MQ\gamma=%
{\displaystyle\sum\limits_{i=1}^{d}}
M_{i}Q\gamma_{i}\in\mathcal{M}_{n\times n}\left(  \mathbb{R}\right)  $.

We denote by $L^{\ast}$ the transpose of the matrix $L$ and $M^{\ast}=\left(
M_{1}^{\ast},...,M_{d}^{\ast}\right)  $.

\section{Formulation of the problem}

Let $\left(  \Omega,\mathcal{F},\left(  \mathcal{F}_{t}\right)  _{t\geq
0},\mathcal{P}\right)  $ be a filtered probability space satisfying the usual
conditions, on which a $d$-dimensional Brownian motion $W=\left(
W_{t}\right)  _{t\geq0}$\ is defined. We assume that $\left(  \mathcal{F}%
_{t}\right)  $ is the $\mathcal{P}$- augmentation of the natural filtration of
$W.$

Let $T$ be a strictly positive real number and $U$ a subset of $\mathbb{R}%
^{m}$.

\subsection{The strict control problem}

\begin{definition}
\textit{An admissible strict control is an }$\mathcal{F}_{t}-$%
\textit{\ adapted process }$v=\left(  v_{t}\right)  $ \textit{with values in
}$U$\textit{\ such that }
\[
\mathbb{E}\left[  \underset{0\leq t\leq T}{\sup}\left\vert v_{t}\right\vert
^{2}\right]  <\infty.
\]

\textit{We denote by }$\mathcal{U}$\textit{\ the set of all admissible strict
controls.}
\end{definition}

For any $v\in\mathcal{U}$, we consider the following controlled BSDE
\begin{equation}
\left\{
\begin{array}
[c]{l}%
dy_{t}^{v}=b\left(  t,y_{t}^{v},z_{t}^{v},v_{t}\right)  dt+z_{t}^{v}dW_{t},\\
y_{T}^{v}=\xi,
\end{array}
\right.
\end{equation}
where $b:\left[  0,T\right]  \times\mathbb{R}^{n}\times\mathcal{M}_{n\times
d}\left(  \mathbb{R}\right)  \times U\longrightarrow\mathbb{R}^{n}$ and $\xi
$\ is an $n$-dimensional $\mathcal{F}_{T}$-measurable random variable such
that%
\[
\mathbb{E}\left\vert \xi\right\vert ^{2}<\infty.
\]

The criteria to be minimized is defined from $\mathcal{U}$ into $\mathbb{R}$
by%
\begin{equation}
J\left(  v\right)  =\mathbb{E}\left[  g\left(  y_{0}^{v}\right)  +%
{\displaystyle\int\nolimits_{0}^{T}}
h\left(  t,y_{t}^{v},z_{t}^{v},v_{t}\right)  dt\right]  ,
\end{equation}
where,%
\begin{align*}
g  &  :\mathbb{R}^{n}\longrightarrow\mathbb{R},\\
h  &  :\left[  0,T\right]  \times\mathbb{R}^{n}\times\mathcal{M}_{n\times
d}\left(  \mathbb{R}\right)  \times U\longrightarrow\mathbb{R}.
\end{align*}

A strict control $u$ is called optimal if it satisfies%
\begin{equation}
J\left(  u\right)  =\inf\limits_{v\in\mathcal{U}}J\left(  v\right)  .
\end{equation}

\ 

We assume that%
\begin{align}
&  b,g\text{ and }h\text{ are continuously differentiable with respect to
}\left(  y,z\right)  \text{.}\\
&  \text{They and all their derivatives with respect to }\left(  y,z\right)
\text{ are continuous in }\left(  y,z,v\right)  .\nonumber\\
&  \text{They are bounded by }C\ \left(  1+\left\vert y\right\vert +\left\vert
z\right\vert +\left\vert v\right\vert \right)  \text{ and their derivatives}%
\nonumber\\
&  \text{with respect to }\left(  y,z\right)  \text{ are continuous and
uniformly bounded.}\nonumber
\end{align}

Under the above hypothesis, for every $v\in U$, equation $\left(  1\right)  $
has a unique strong solution and the functional cost $J$ is well defined from
$\mathcal{U}$ into $\mathbb{R}$.

\subsection{The relaxed model}

The idea for relaxed the strict control problem defined above is to embed the
set $U$ of strict controls into a wider class which gives a more suitable
topological structure. In the relaxed model, the $U$-valued process $v$ is
replaced by a $\mathbb{P}\left(  U\right)  $-valued process $q$, where
$\mathbb{P}\left(  U\right)  $ denotes the space of probability measure on $U$
equipped with the topology of stable convergence.

\begin{definition}
A relaxed control $\left(  q_{t}\right)  _{t}$ is a $\mathbb{P}\left(
U\right)  $-valued process, progressively measurable with respect to $\left(
\mathcal{F}_{t}\right)  _{t}$\ and such that for each $t$, $1_{]0,t]}.q$\ is
$\mathcal{F}_{t}$-measurable.

\textit{We denote by }$\mathcal{R}$\textit{\ the set of all relaxed controls.}
\end{definition}

\begin{remark}
The set of strict controls is embedded into the set\ of relaxed controls by
the mapping
\[
f:v\mathbb{\longmapsto}f_{v}\left(  dt,da\right)  =dt\delta_{v_{t}}(da),
\]
where $\delta_{v}$ is the atomic measure concentrated at a single point $v$.
\end{remark}

For more details on relaxed controls, see $\left[  1\right]  ,\left[
3\right]  ,\left[  4\right]  ,\left[  5\right]  ,\left[  10\right]  ,\left[
14\right]  ,\left[  21\right]  ,\left[  23\right]  ,\left[  24\right]  .$

\ 

For any $q\in\mathcal{R}$, we consider the following relaxed BSDE%
\begin{equation}
\left\{
\begin{array}
[c]{l}%
dy_{t}^{q}=\int_{U}b\left(  t,y_{t}^{q},z_{t}^{q},a\right)  q_{t}\left(
da\right)  dt+z_{t}^{q}dW_{t},\\
y_{T}^{q}=\xi.
\end{array}
\right.
\end{equation}

The expected cost to be minimized, in the relaxed model, is defined from
$\mathcal{R}$ into $\mathbb{R}$ by%
\begin{equation}
J\left(  q\right)  =\mathbb{E}\left[  g\left(  y_{0}^{q}\right)  +%
{\displaystyle\int\nolimits_{0}^{T}}
\int_{U}h\left(  t,y_{t}^{q},z_{t}^{q},a\right)  q_{t}\left(  da\right)
dt\right]  .
\end{equation}

A relaxed control $\mu$ is called optimal if it solves%
\begin{equation}
J\left(  \mu\right)  =\inf\limits_{q\in\mathcal{R}}J\left(  q\right)  .
\end{equation}

\begin{remark}
If we put
\begin{align*}
\overline{b}\left(  t,y_{t}^{q},z_{t}^{q},q_{t}\right)   &  =\int_{U}b\left(
t,y_{t}^{q},z_{t}^{q},a\right)  q_{t}\left(  da\right)  ,\\
\overline{h}\left(  t,y_{t}^{q},z_{t}^{q},q_{t}\right)   &  =\int_{U}h\left(
t,y_{t}^{q},z_{t}^{q},a\right)  q_{t}\left(  da\right)  .
\end{align*}

Then, equation $\left(  5\right)  $ becomes
\begin{equation}
\left\{
\begin{array}
[c]{l}%
dy_{t}^{q}=\overline{b}\left(  t,y_{t}^{q},z_{t}^{q},q_{t}\right)
dt+z_{t}^{q}dW_{t},\\
y_{T}^{q}=\xi.
\end{array}
\right.  \tag{5$^\prime$}%
\end{equation}

With a functional cost given by%
\[
J\left(  q\right)  =\mathbb{E}\left[  g\left(  y_{0}^{q}\right)  +%
{\displaystyle\int\nolimits_{0}^{T}}
\overline{h}\left(  t,y_{t}^{q},z_{t}^{q},q_{t}\right)  dt\right]  .
\]

Hence, by introducing relaxed controls, we have replaced $U$ by a larger space
$\mathbb{P}\left(  U\right)  $. We have gained the advantage that
$\mathbb{P}\left(  U\right)  $ is both compact and convex. Furthermore, the
new coefficients of equation $\left(  5^{\prime}\right)  $ and the running
cost are linear with respect to the relaxed control variable.
\end{remark}

\begin{remark}
The coefficient $\overline{b}$ (defined in the above remark) check
respectively the same assumptions as $b$. Then, under assumptions $\left(
4\right)  $, $\overline{b}$ is uniformly Lipschitz and with linear growth.
Then by classical results on BSDEs, for every $q\in\mathcal{R}$ equation
$\left(  5^{\prime}\right)  $ admits a unique strong solution. Consequently,
for every $q\in\mathcal{R}$ equation $\left(  5\right)  $ has a unique strong solution.

On the other hand, It is easy to see that $\overline{h}$ checks the same
assumptions as $h$. Then, the functional cost $J$ is well defined from
$\mathcal{R}$ into $\mathbb{R}$.
\end{remark}

\begin{remark}
If $q_{t}=\delta_{v_{t}}$ is an atomic measure concentrated at a single point
$v_{t}\in U$, then for each $t\in\left[  0,T\right]  $ we have
\begin{align*}
\int_{U}b\left(  t,y_{t}^{q},z_{t}^{q},a\right)  q_{t}\left(  da\right)   &
=\int_{U}b\left(  t,y_{t}^{q},z_{t}^{q},a\right)  \delta_{v_{t}}\left(
da\right)  =b\left(  t,y_{t}^{q},z_{t}^{q},v_{t}\right)  ,\\
\int_{U}h\left(  t,y_{t}^{q},z_{t}^{q},a\right)  q_{t}\left(  da\right)   &
=\int_{U}h\left(  t,y_{t}^{q},z_{t}^{q},a\right)  \delta_{v_{t}}\left(
da\right)  =h\left(  t,y_{t}^{q},z_{t}^{q},v_{t}\right)  .
\end{align*}

In this case $\left(  y^{q},z^{q}\right)  =\left(  y^{v},z^{v}\right)  $,
$J\left(  q\right)  =J\left(  v\right)  $ and we get a strict control problem.
So the problem of strict controls $\left\{  \left(  1\right)  ,\left(
2\right)  ,\left(  3\right)  \right\}  $ is a particular case of relaxed
control problem $\left\{  \left(  5\right)  ,\left(  6\right)  ,\left(
7\right)  \right\}  $.
\end{remark}

\section{Necessary and sufficient optimality conditions for relaxed controls}

In this section, we study the problem $\left\{  \left(  5\right)  ,\left(
6\right)  ,\left(  7\right)  \right\}  $ and we establish necessary as well as
sufficient conditions of optimality for relaxed controls.

\subsection{Preliminary results}

Since the set $\mathcal{R}$ is convex, then the classical way to derive
necessary optimality conditions for relaxed controls is to use the convex
perturbation method. More precisely, let $\mu$ be an optimal relaxed control
and $\left(  y_{t}^{\mu},z_{t}^{\mu}\right)  $ the solution of $\left(
5\right)  $ controlled by $\mu$. Then, we can define a perturbed relaxed
control as follows%
\begin{equation}
\mu_{t}^{\theta}=\mu_{t}+\theta\left(  q_{t}-\mu_{t}\right)  ,
\end{equation}
where, $\theta>0$ is sufficiently small and $q$ is an arbitrary element of
$\mathcal{R}$.

Denote by $\left(  y_{t}^{\theta},z_{t}^{\theta}\right)  $ the solution of
$\left(  5\right)  $ associated with $\mu^{\theta}$.

From optimality of $\mu$, the variational inequality will be derived from the
fact that
\begin{equation}
0\leq J\left(  \mu^{\theta}\right)  -J\left(  \mu\right)  .
\end{equation}

For this end, we need the following classical lemmas.

\begin{lemma}
\textit{Under assumptions }$\left(  4\right)  $, we have%
\begin{align}
\underset{\theta\rightarrow0}{\lim}\left[  \underset{0\leq t\leq T}{\sup
}\mathbb{E}\left\vert y_{t}^{\theta}-y_{t}^{\mu}\right\vert ^{2}\right]   &
=0,\\
\underset{\theta\rightarrow0}{\lim}\mathbb{E}%
{\displaystyle\int\nolimits_{0}^{T}}
\left\vert z_{t}^{\theta}-z_{t}^{\mu}\right\vert ^{2}dt  &  =0.
\end{align}

\end{lemma}

\begin{proof}
Applying It\^{o}'s formula to $\left(  y_{t}^{\theta}-y_{t}^{\mu}\right)
^{2}$, we have%
\begin{align*}
&  \mathbb{E}\left\vert y_{t}^{\theta}-y_{t}^{\mu}\right\vert ^{2}+\mathbb{E}%
{\displaystyle\int\nolimits_{t}^{T}}
\left\vert z_{s}^{\theta}-z_{s}^{\mu}\right\vert ^{2}ds\\
&  =2\mathbb{E}%
{\displaystyle\int\nolimits_{t}^{T}}
\left\vert \left(  y_{s}^{\theta}-y_{s}^{\mu}\right)  \left[
{\displaystyle\int\nolimits_{U}}
b\left(  s,y_{s}^{\theta},z_{s}^{\theta},a\right)  \mu_{s}^{\theta}\left(
da\right)  -%
{\displaystyle\int\nolimits_{U}}
b\left(  s,y_{s}^{\mu},z_{s}^{\mu},a\right)  \mu_{s}\left(  da\right)
\right]  \right\vert \,ds.
\end{align*}

From the Young formula, for every $\varepsilon>0$, we have%
\begin{align*}
&  \mathbb{E}\left\vert y_{t}^{\theta}-y_{t}^{\mu}\right\vert ^{2}+\mathbb{E}%
{\displaystyle\int\nolimits_{t}^{T}}
\left\vert z_{s}^{\theta}-z_{s}^{\mu}\right\vert ^{2}ds\\
&  \leq%
\genfrac{.}{.}{}{0}{1}{\varepsilon}%
\mathbb{E}%
{\displaystyle\int\nolimits_{t}^{T}}
\left\vert y_{s}^{\theta}-y_{s}^{\mu}\right\vert ^{2}ds\\
&  +\varepsilon\mathbb{E}%
{\displaystyle\int\nolimits_{t}^{T}}
\left\vert
{\displaystyle\int\nolimits_{U}}
b\left(  s,y_{s}^{\theta},z_{s}^{\theta},a\right)  \mu_{s}^{\theta}\left(
da\right)  -%
{\displaystyle\int\nolimits_{U}}
b\left(  s,y_{s}^{\mu},z_{s}^{\mu},a\right)  \mu_{s}\left(  da\right)
\right\vert ^{2}ds.
\end{align*}

Then,%
\begin{align*}
&  \mathbb{E}\left\vert y_{t}^{\theta}-y_{t}^{\mu}\right\vert ^{2}+\mathbb{E}%
{\displaystyle\int\nolimits_{t}^{T}}
\left\vert z_{s}^{\theta}-z_{s}^{\mu}\right\vert ^{2}ds\\
&  \leq%
\genfrac{.}{.}{}{0}{1}{\varepsilon}%
\mathbb{E}%
{\displaystyle\int\nolimits_{t}^{T}}
\left\vert y_{s}^{\theta}-y_{s}^{\mu}\right\vert ^{2}ds\\
&  +C\varepsilon\mathbb{E}%
{\displaystyle\int\nolimits_{t}^{T}}
\left\vert
{\displaystyle\int\nolimits_{U}}
b\left(  s,y_{s}^{\theta},z_{s}^{\theta},a\right)  \mu_{s}^{\theta}\left(
da\right)  -%
{\displaystyle\int\nolimits_{U}}
b\left(  s,y_{s}^{\theta},z_{s}^{\theta},a\right)  \mu_{s}\left(  da\right)
\right\vert ^{2}ds\\
&  +C\varepsilon\mathbb{E}%
{\displaystyle\int\nolimits_{t}^{T}}
\left\vert
{\displaystyle\int\nolimits_{U}}
b\left(  s,y_{s}^{\theta},z_{s}^{\theta},a\right)  \mu_{s}\left(  da\right)  -%
{\displaystyle\int\nolimits_{U}}
b\left(  s,y_{s}^{\mu},z_{s}^{\theta},a\right)  \mu_{s}\left(  da\right)
\right\vert ^{2}ds\\
&  +C\varepsilon\mathbb{E}%
{\displaystyle\int\nolimits_{t}^{T}}
\left\vert
{\displaystyle\int\nolimits_{U}}
b\left(  s,y_{s}^{\mu},z_{s}^{\theta},a\right)  \mu_{s}\left(  da\right)  -%
{\displaystyle\int\nolimits_{U}}
b\left(  s,y_{s}^{\mu},z_{s}^{\mu},a\right)  \mu_{s}\left(  da\right)
\right\vert ds.
\end{align*}

By the definition of $\mu_{t}^{\theta}$, we have%
\begin{align*}
&  \mathbb{E}\left\vert y_{t}^{\theta}-y_{t}^{\mu}\right\vert ^{2}+\mathbb{E}%
{\displaystyle\int\nolimits_{t}^{T}}
\left\vert z_{s}^{\theta}-z_{s}^{\mu}\right\vert ^{2}ds\\
&  \leq%
\genfrac{.}{.}{}{0}{1}{\varepsilon}%
\mathbb{E}%
{\displaystyle\int\nolimits_{t}^{T}}
\left\vert y_{s}^{\theta}-y_{s}^{\mu}\right\vert ^{2}ds\\
&  +C\varepsilon\theta^{2}\mathbb{E}%
{\displaystyle\int\nolimits_{t}^{T}}
\left\vert
{\displaystyle\int\nolimits_{U}}
b\left(  s,y_{s}^{\theta},z_{s}^{\theta},a\right)  q_{s}\left(  da\right)  -%
{\displaystyle\int\nolimits_{U}}
b\left(  s,y_{s}^{\theta},z_{s}^{\theta},a\right)  \mu_{s}\left(  da\right)
\right\vert ^{2}ds\\
&  +C\varepsilon\mathbb{E}%
{\displaystyle\int\nolimits_{t}^{T}}
\left\vert
{\displaystyle\int\nolimits_{U}}
b\left(  s,y_{s}^{\theta},z_{s}^{\theta},a\right)  \mu_{s}\left(  da\right)  -%
{\displaystyle\int\nolimits_{U}}
b\left(  s,y_{s}^{\mu},z_{s}^{\theta},a\right)  \mu_{s}\left(  da\right)
\right\vert ^{2}ds\\
&  +C\varepsilon\mathbb{E}%
{\displaystyle\int\nolimits_{t}^{T}}
\left\vert
{\displaystyle\int\nolimits_{U}}
b\left(  s,y_{s}^{\mu},z_{s}^{\theta},a\right)  \mu_{s}\left(  da\right)  -%
{\displaystyle\int\nolimits_{U}}
b\left(  s,y_{s}^{\mu},z_{s}^{\mu},a\right)  \mu_{s}\left(  da\right)
\right\vert ^{2}ds.
\end{align*}

Since $b$\ is uniformly Lipschitz with respect to $y,z$, then%
\begin{align*}
\mathbb{E}\left\vert y_{t}^{\theta}-y_{t}^{\mu}\right\vert ^{2}+\mathbb{E}%
{\displaystyle\int\nolimits_{t}^{T}}
\left\vert z_{s}^{\theta}-z_{s}^{\mu}\right\vert ^{2}ds  &  \leq\left(
\genfrac{.}{.}{}{0}{1}{\varepsilon}%
\mathbb{+}C\ \varepsilon\right)  \mathbb{E}%
{\displaystyle\int\nolimits_{t}^{T}}
\left\vert y_{s}^{\theta}-y_{s}^{\mu}\right\vert ^{2}ds\\
&  +C\varepsilon\mathbb{E}%
{\displaystyle\int\nolimits_{t}^{T}}
\left\vert z_{s}^{\theta}-z_{s}^{\mu}\right\vert ^{2}ds+C\varepsilon\theta
^{2}.
\end{align*}

Choose $\varepsilon=%
\genfrac{.}{.}{}{0}{1}{2C}%
$, then we have%
\[
\mathbb{E}\left\vert y_{t}^{\theta}-y_{t}^{\mu}\right\vert ^{2}+%
\genfrac{.}{.}{}{0}{1}{2}%
\mathbb{E}%
{\displaystyle\int\nolimits_{t}^{T}}
\left\vert z_{s}^{\theta}-z_{s}^{\mu}\right\vert ^{2}ds\leq\left(  2C+%
\genfrac{.}{.}{}{0}{1}{2}%
\right)  \mathbb{E}%
{\displaystyle\int\nolimits_{t}^{T}}
\left\vert y_{s}^{\theta}-y_{s}^{\mu}\right\vert ^{2}ds+C\varepsilon\theta
^{2}.
\]

From the above inequality, we derive two inequalities%
\begin{equation}
\mathbb{E}\left\vert y_{t}^{\theta}-y_{t}^{\mu}\right\vert ^{2}\leq\left(  2C+%
\genfrac{.}{.}{}{0}{1}{2}%
\right)  \mathbb{E}%
{\displaystyle\int\nolimits_{t}^{T}}
\left\vert y_{s}^{\theta}-y_{s}^{\mu}\right\vert ^{2}ds+C\varepsilon\theta
^{2},
\end{equation}%
\begin{equation}
\mathbb{E}%
{\displaystyle\int\nolimits_{t}^{T}}
\left\vert z_{s}^{\theta}-z_{s}^{\mu}\right\vert ^{2}ds\leq\left(
4C+1\right)  \mathbb{E}%
{\displaystyle\int\nolimits_{t}^{T}}
\left\vert y_{s}^{\theta}-y_{s}^{\mu}\right\vert ^{2}ds+2C\varepsilon
\theta^{2}.
\end{equation}

By using $\left(  12\right)  $, Gronwall's lemma and Bukholder-Davis-Gundy
inequality, we obtain $\left(  10\right)  $.\ Finally, $\left(  11\right)
$\ is derived from $\left(  10\right)  $\ and $\left(  13\right)  $.
\end{proof}

\begin{lemma}
\textit{Let }$\widetilde{y}_{t}$ be \textit{the solution of the following
linear equation (called variational equation)}\
\begin{equation}
\left\{
\begin{array}
[c]{ll}%
d\widetilde{y}_{t}= &
{\displaystyle\int\nolimits_{U}}
\left[  b_{y}\left(  t,y_{t}^{\mu},z_{t}^{\mu},a\right)  \widetilde{y}%
_{t}+b_{z}\left(  t,y_{t}^{\mu},z_{t}^{\mu},a\right)  \widetilde{z}%
_{t}\right]  \mu_{t}\left(  da\right)  dt\\
& +\left[
{\displaystyle\int\nolimits_{U}}
b\left(  t,y_{t}^{\mu},z_{t}^{\mu},a\right)  q_{t}\left(  da\right)  -%
{\displaystyle\int\nolimits_{U}}
b\left(  t,y_{t}^{\mu},z_{t}^{\mu},a\right)  \mu_{t}\left(  da\right)
\right]  dt+\widetilde{z}_{t}dW_{t},\\
\widetilde{y}_{T}= & 0.
\end{array}
\right.
\end{equation}

T\textit{hen, the following estimations hold}%
\begin{align}
\underset{\theta\rightarrow0}{\lim}\mathbb{E}\left\vert
\genfrac{.}{.}{}{0}{y_{t}^{\theta}-y_{t}^{\mu}}{\theta}%
-\widetilde{y}_{t}\right\vert ^{2}  &  =0,\\
\underset{\theta\rightarrow0}{\lim}\mathbb{E}%
{\displaystyle\int\nolimits_{0}^{T}}
\left\vert
\genfrac{.}{.}{}{0}{z_{t}^{\theta}-z_{t}^{\mu}}{\theta}%
-\widetilde{z}_{t}\right\vert ^{2}dt  &  =0.
\end{align}

\end{lemma}

\begin{proof}
For simplicity, we put%
\begin{align}
Y_{t}^{\theta}  &  =%
\genfrac{.}{.}{}{0}{y_{t}^{\theta}-y_{t}^{\mu}}{\theta}%
-\widetilde{y}_{t},\\
Z_{t}^{\theta}  &  =%
\genfrac{.}{.}{}{0}{z_{t}^{\theta}-z_{t}^{\mu}}{\theta}%
-\widetilde{z}_{t}.
\end{align}%
\begin{equation}
\Lambda_{t}^{\theta}\left(  a\right)  =\left(  t,y_{t}^{\mu}+\lambda
\theta\left(  Y_{t}^{\theta}+\widetilde{y}_{t}\right)  ,z_{t}^{\mu}%
+\lambda\theta\left(  Z_{t}^{\theta}+\widetilde{z}_{t}\right)  ,a\right)  .
\end{equation}

By $\left(  17\right)  $ and\ $\left(  18\right)  $ we have the following BSDE%
\[
\left\{
\begin{array}
[c]{l}%
dY_{t}^{\theta}=\left(  F_{t}^{y}Y_{t}^{\theta}dt+F_{t}^{y}Z_{t}^{\theta
}-\gamma_{t}^{\theta}\right)  dt+Z_{t}^{\theta}dW_{t},\\
Y_{T}^{\theta}=0,
\end{array}
\right.
\]
where,%
\begin{align*}
F_{t}^{y}  &  =-\int_{0}^{1}\int_{U}b_{y}\left(  \Lambda_{t}^{\theta}\left(
a\right)  \right)  \mu_{t}\left(  da\right)  d\lambda,\\
F_{t}^{z}  &  =-\int_{0}^{1}\int_{U}b_{z}\left(  \Lambda_{t}^{\theta}\left(
a\right)  \right)  \mu_{t}\left(  da\right)  d\lambda,
\end{align*}
and $\gamma_{t}^{\theta}$ is given by%
\begin{align*}
&  \gamma_{t}^{\theta}=\int_{t}^{T}\int_{U}\left[  b_{y}\left(  \Lambda
_{s}^{\theta}\left(  a\right)  \right)  \left(  y_{s}^{\theta}-y_{s}^{\mu
}\right)  +b_{z}\left(  \Lambda_{s}^{\theta}\left(  a\right)  \right)  \left(
z_{s}^{\theta}-z_{s}^{\mu}\right)  \right]  q_{s}\left(  da\right)  ds\\
&  -\int_{t}^{T}\int_{U}\left[  b_{y}\left(  \Lambda_{s}^{\theta}\left(
a\right)  \right)  \left(  y_{s}^{\theta}-y_{s}^{\mu}\right)  +b_{z}\left(
\Lambda_{s}^{\theta}\left(  a\right)  \right)  \left(  z_{s}^{\theta}%
-z_{s}^{\mu}\right)  \right]  \mu_{s}\left(  da\right)  ds.
\end{align*}

Since $b_{y}$ and $b_{z}$ are bounded, then%
\[
\mathbb{E}\left\vert \gamma_{t}^{\theta}\right\vert ^{2}\leq C\mathbb{E}%
\int_{t}^{T}\left\vert y_{s}^{\theta}-y_{s}^{\mu}\right\vert ^{2}%
ds+C\mathbb{E}\int_{t}^{T}\left\vert z_{s}^{\theta}-z_{s}^{\mu}\right\vert
^{2}ds
\]

By $\left(  10\right)  $ and $\left(  11\right)  $, we get%
\begin{equation}
\underset{\theta\rightarrow0}{\lim}\mathbb{E}\left\vert \gamma_{t}^{\theta
}\right\vert ^{2}=0.
\end{equation}

Applying It\^{o}'s formula to $\left(  Y_{t}^{\theta}\right)  ^{2}$, we get%
\[
\mathbb{E}\left\vert Y_{t}^{\theta}\right\vert ^{2}+\mathbb{E}%
{\displaystyle\int\nolimits_{t}^{T}}
\left\vert Z_{s}^{\theta}\right\vert ^{2}ds=\mathbb{E}\left\vert Y_{T}%
^{\theta}\right\vert ^{2}+2\mathbb{E}%
{\displaystyle\int\nolimits_{t}^{T}}
\left\vert Y_{s}^{\theta}\left(  F_{s}^{y}Y_{s}^{\theta}+F_{s}^{z}%
Z_{s}^{\theta}-\gamma_{s}^{\theta}\right)  \right\vert ds.
\]

By using the Young formula, for every $\varepsilon>0$, we have%
\begin{align*}
\mathbb{E}\left\vert Y_{t}^{\theta}\right\vert ^{2}+\mathbb{E}%
{\displaystyle\int\nolimits_{t}^{T}}
\left\vert Z_{s}^{\theta}\right\vert ^{2}ds  &  \leq\mathbb{E}\left\vert
Y_{T}^{\theta}\right\vert ^{2}+%
\genfrac{.}{.}{}{0}{1}{\varepsilon}%
\mathbb{E}%
{\displaystyle\int\nolimits_{t}^{T}}
\left\vert Y_{s}^{\theta}\right\vert ^{2}ds+\varepsilon\mathbb{E}%
{\displaystyle\int\nolimits_{t}^{T}}
\left\vert \left(  F_{s}^{y}Y_{s}^{\theta}+F_{s}^{z}Z_{s}^{\theta}-\gamma
_{s}^{\theta}\right)  \right\vert ^{2}ds\\
&  \leq\mathbb{E}\left\vert Y_{T}^{\theta}\right\vert ^{2}+%
\genfrac{.}{.}{}{0}{1}{\varepsilon}%
\mathbb{E}%
{\displaystyle\int\nolimits_{t}^{T}}
\left\vert Y_{s}^{\theta}\right\vert ^{2}ds+C\varepsilon\mathbb{E}%
{\displaystyle\int\nolimits_{t}^{T}}
\left\vert F_{s}^{y}Y_{s}^{\theta}\right\vert ^{2}ds\\
&  +C\varepsilon\mathbb{E}%
{\displaystyle\int\nolimits_{t}^{T}}
\left\vert F_{s}^{z}Z_{s}^{\theta}\right\vert ^{2}ds+C\varepsilon\mathbb{E}%
{\displaystyle\int\nolimits_{t}^{T}}
\left\vert \gamma_{s}^{\theta}\right\vert ^{2}ds.
\end{align*}

Since $F_{t}^{y}$ and $F_{t}^{z}$ are bounded, then%
\[
\mathbb{E}\left\vert Y_{t}^{\theta}\right\vert ^{2}+\mathbb{E}%
{\displaystyle\int\nolimits_{t}^{T}}
\left\vert Z_{s}^{\theta}\right\vert ^{2}ds\leq\left(
\genfrac{.}{.}{}{0}{1}{\varepsilon}%
+C\ \varepsilon\right)  \mathbb{E}%
{\displaystyle\int\nolimits_{t}^{T}}
\left\vert Y_{s}^{\theta}\right\vert ^{2}ds+C\ \varepsilon\mathbb{E}%
{\displaystyle\int\nolimits_{t}^{T}}
\left\vert Z_{s}^{\theta}\right\vert ^{2}ds+C\ \varepsilon\mathbb{E}%
{\displaystyle\int\nolimits_{t}^{T}}
\left\vert \gamma_{s}^{\theta}\right\vert ^{2}ds,
\]

Choose $\varepsilon=%
\genfrac{.}{.}{}{0}{1}{2C}%
$, then we have%
\[
\mathbb{E}\left\vert Y_{t}^{\theta}\right\vert ^{2}+%
\genfrac{.}{.}{}{0}{1}{2}%
\mathbb{E}%
{\displaystyle\int\nolimits_{t}^{T}}
\left\vert Z_{s}^{\theta}\right\vert ^{2}ds\leq\left(  2C+%
\genfrac{.}{.}{}{0}{1}{2}%
\right)  \mathbb{E}%
{\displaystyle\int\nolimits_{t}^{T}}
\left\vert Y_{s}^{\theta}\right\vert ^{2}ds+C\ \varepsilon\mathbb{E}%
{\displaystyle\int\nolimits_{t}^{T}}
\left\vert \gamma_{s}^{\theta}\right\vert ^{2}ds.
\]

From the above inequality, we deduce two inequalities%
\begin{equation}
\mathbb{E}\left\vert Y_{t}^{\theta}\right\vert ^{2}\leq\left(  2C+%
\genfrac{.}{.}{}{0}{1}{2}%
\right)  \mathbb{E}%
{\displaystyle\int\nolimits_{t}^{T}}
\left\vert Y_{s}^{\theta}\right\vert ^{2}ds+C\ \varepsilon\mathbb{E}%
{\displaystyle\int\nolimits_{t}^{T}}
\left\vert \gamma_{s}^{\theta}\right\vert ^{2}ds,
\end{equation}%
\begin{equation}
\mathbb{E}%
{\displaystyle\int\nolimits_{t}^{T}}
\left\vert Z_{s}^{\theta}\right\vert ^{2}ds\leq\left(  4C+1\right)  \mathbb{E}%
{\displaystyle\int\nolimits_{t}^{T}}
\left\vert Y_{s}^{\theta}\right\vert ^{2}ds+2C\ \varepsilon\mathbb{E}%
{\displaystyle\int\nolimits_{t}^{T}}
\left\vert \gamma_{s}^{\theta}\right\vert ^{2}ds.
\end{equation}

By using $\left(  20\right)  ,\ \left(  21\right)  $ and Gronwall's lemma, we
obtain $\left(  15\right)  $. Finally $\left(  16\right)  $ is derived from
$\left(  15\right)  ,\ \left(  20\right)  $ and $\left(  21\right)  $.
\end{proof}

\begin{lemma}
\textit{Let }$\mu$\textit{\ be an optimal control minimizing the functional
}$J$\textit{\ over }$\mathcal{R}$\textit{\ and }$\left(  y_{t}^{\mu}%
,z_{t}^{\mu}\right)  $\textit{\ the solution of }$\left(  5\right)
$\textit{\ associated with }$\mu$\textit{. Then for any }$q\in\mathcal{R}%
$\textit{, we have}
\begin{align}
0  &  \leq\mathbb{E}\left[  g_{y}\left(  y_{0}^{\mu}\right)  \widetilde{y}%
_{0}\right] \\
&  +\mathbb{E}%
{\displaystyle\int\nolimits_{0}^{T}}
\left[
{\displaystyle\int\nolimits_{U}}
h\left(  t,y_{t}^{\mu},z_{t}^{\mu},a\right)  q_{t}\left(  da\right)  -%
{\displaystyle\int\nolimits_{U}}
h\left(  t,y_{t}^{\mu},z_{t}^{\mu},a\right)  \mu_{t}\left(  da\right)
\right]  dt\nonumber\\
&  +\mathbb{E}%
{\displaystyle\int\nolimits_{0}^{T}}
{\displaystyle\int\nolimits_{U}}
h_{y}\left(  t,y_{t}^{\mu},z_{t}^{\mu},a\right)  \mu_{t}\left(  da\right)
\widetilde{y}_{t}dt+\mathbb{E}%
{\displaystyle\int\nolimits_{0}^{T}}
{\displaystyle\int\nolimits_{U}}
h_{z}\left(  t,y_{t}^{\mu},z_{t}^{\mu},a\right)  \mu_{t}\left(  da\right)
\widetilde{z}_{t}dt.\nonumber
\end{align}

\end{lemma}

\begin{proof}
Let $\mu$ be an optimal relaxed control minimizing the cost $J$ over
$\mathcal{R}$, then from $\left(  9\right)  $ we have
\begin{align*}
0  &  \leq\mathbb{E}\left[  g\left(  y_{0}^{\theta}\right)  -g\left(
y_{0}^{\mu}\right)  \right] \\
&  +\mathbb{E}%
{\displaystyle\int\nolimits_{0}^{T}}
\left[
{\displaystyle\int\nolimits_{U}}
h\left(  t,y_{t}^{\theta},z_{t}^{\theta},a\right)  \mu_{t}^{\theta}\left(
da\right)  -%
{\displaystyle\int\nolimits_{U}}
h\left(  t,y_{t}^{\mu},z_{t}^{\mu},a\right)  \mu_{t}\left(  da\right)
\right]  dt\\
&  =\mathbb{E}\left[  g\left(  y_{0}^{\theta}\right)  -g\left(  y_{0}^{\mu
}\right)  \right] \\
&  +\mathbb{E}%
{\displaystyle\int\nolimits_{0}^{T}}
\left[
{\displaystyle\int\nolimits_{U}}
h\left(  t,y_{t}^{\theta},z_{t}^{\theta},a\right)  \mu_{t}^{\theta}\left(
da\right)  -%
{\displaystyle\int\nolimits_{U}}
h\left(  t,y_{t}^{\theta},z_{t}^{\theta},a\right)  \mu_{t}\left(  da\right)
\right]  dt\\
&  +\mathbb{E}%
{\displaystyle\int\nolimits_{0}^{T}}
\left[
{\displaystyle\int\nolimits_{U}}
h\left(  t,y_{t}^{\theta},z_{t}^{\theta},a\right)  \mu_{t}\left(  da\right)  -%
{\displaystyle\int\nolimits_{U}}
h\left(  t,y_{t}^{\mu},z_{t}^{\mu},a\right)  \mu_{t}\left(  da\right)
\right]  dt.
\end{align*}

From the definition of $\mu^{\theta}$, we get%
\begin{align*}
0  &  \leq\mathbb{E}\left[  g\left(  y_{0}^{\theta}\right)  -g\left(
y_{0}^{\mu}\right)  \right] \\
&  +\theta\mathbb{E}%
{\displaystyle\int\nolimits_{0}^{T}}
\left[
{\displaystyle\int\nolimits_{U}}
h\left(  t,y_{t}^{\theta},z_{t}^{\theta},a\right)  q_{t}\left(  da\right)  -%
{\displaystyle\int\nolimits_{U}}
h\left(  t,y_{t}^{\theta},z_{t}^{\theta},a\right)  \mu_{t}\left(  da\right)
\right]  dt\\
&  +\mathbb{E}%
{\displaystyle\int\nolimits_{0}^{T}}
{\displaystyle\int\nolimits_{U}}
\left[  h\left(  t,y_{t}^{\theta},z_{t}^{\theta},a\right)  -h\left(
t,y_{t}^{\mu},z_{t}^{\mu},a\right)  \right]  \mu_{t}\left(  da\right)  dt.
\end{align*}

Then,%
\begin{align}
0  &  \leq\mathbb{E}%
{\displaystyle\int\nolimits_{0}^{1}}
\left[  g_{y}\left(  y_{0}^{\mu}+\lambda\theta\left(  \widetilde{y}_{0}%
+Y_{0}^{\theta}\right)  \right)  \widetilde{y}_{0}\right]  d\lambda\\
&  +\mathbb{E}%
{\displaystyle\int\nolimits_{0}^{T}}
{\displaystyle\int\nolimits_{0}^{1}}
{\displaystyle\int\nolimits_{U}}
\left[  h_{y}\left(  \Lambda_{t}^{\theta}\left(  a\right)  \right)
\widetilde{y}_{t}+h_{z}\left(  \Lambda_{t}^{\theta}\left(  a\right)  \right)
\widetilde{z}_{t}\right]  \mu_{t}\left(  da\right)  d\lambda dt\nonumber\\
&  +\mathbb{E}%
{\displaystyle\int\nolimits_{0}^{T}}
\left[
{\displaystyle\int\nolimits_{U}}
h\left(  t,y_{t}^{\mu},z_{t}^{\mu},a\right)  q_{t}\left(  da\right)  -%
{\displaystyle\int\nolimits_{U}}
h\left(  t,y_{t}^{\mu},z_{t}^{\mu},a\right)  \mu_{t}\left(  da\right)
\right]  dt+\rho_{t}^{\theta},\nonumber
\end{align}
where $\rho_{t}^{\theta}$ is given by%
\begin{align*}
&  \rho_{t}^{\theta}=\mathbb{E}%
{\displaystyle\int\nolimits_{0}^{1}}
\left[  g_{y}\left(  y_{0}^{\mu}+\lambda\theta\left(  \widetilde{y}_{0}%
+Y_{0}^{\theta}\right)  \right)  Y_{0}^{\theta}\right]  d\lambda\\
&  +\mathbb{E}%
{\displaystyle\int\nolimits_{0}^{T}}
{\displaystyle\int\nolimits_{0}^{1}}
{\displaystyle\int\nolimits_{U}}
\left[  h_{y}\left(  \Lambda_{t}^{\theta}\left(  a\right)  \right)  \left(
y_{t}^{\theta}-y_{t}^{\mu}\right)  +h_{z}\left(  \Lambda_{t}^{\theta}\left(
a\right)  \right)  \left(  z_{t}^{\theta}-z_{t}^{\mu}\right)  \right]  \mu
_{t}\left(  da\right)  d\lambda dt\\
&  +\mathbb{E}%
{\displaystyle\int\nolimits_{0}^{T}}
{\displaystyle\int\nolimits_{0}^{1}}
{\displaystyle\int\nolimits_{U}}
\left[  h_{y}\left(  \Lambda_{t}^{\theta}\left(  a\right)  \right)
Y_{t}^{\theta}+h_{z}\left(  \Lambda_{t}^{\theta}\left(  a\right)  \right)
Z_{t}^{\theta}\right]  \mu_{t}\left(  da\right)  d\lambda dt.
\end{align*}

Since the derivatives $g_{y},h_{y},h_{z}$ are bounded, then by using the
Cauchy-Schwartz inequality, we have%
\begin{align*}
\rho_{t}^{\theta}  &  =\left(  \mathbb{E}\left\vert Y_{0}^{\theta}\right\vert
^{2}\right)  ^{1/2}+\left(  \mathbb{E}%
{\displaystyle\int\nolimits_{0}^{T}}
\left\vert y_{t}^{\theta}-y_{t}^{\mu}\right\vert ^{2}dt\right)  ^{1/2}+\left(
\mathbb{E}%
{\displaystyle\int\nolimits_{0}^{T}}
\left\vert z_{t}^{\theta}-z_{t}^{\mu}\right\vert ^{2}dt\right)  ^{1/2}\\
&  +\left(  \mathbb{E}%
{\displaystyle\int\nolimits_{0}^{T}}
\left\vert Y_{t}^{\theta}\right\vert ^{2}dt\right)  ^{1/2}+\left(  \mathbb{E}%
{\displaystyle\int\nolimits_{0}^{T}}
\left\vert Z_{t}^{\theta}\right\vert ^{2}dt\right)  ^{1/2}.
\end{align*}

By $\left(  10\right)  ,\ \left(  11\right)  ,\ \left(  15\right)  $
and$\ \left(  16\right)  $, we get
\[
\underset{\theta\rightarrow0}{\lim}\rho_{t}^{\theta}=0.
\]

Finally, by letting $\theta$ go to $0$ in $\left(  24\right)  $, the proof is completed.
\end{proof}

\subsection{Necessary optimality conditions for relaxed controls}

Starting from the variational inequality $\left(  23\right)  $, we can now
state necessary optimality conditions for the relaxed control problem
$\left\{  \left(  5\right)  ,\left(  6\right)  ,\left(  7\right)  \right\}  $
in the global form.

\begin{theorem}
(Necessary optimality conditions for relaxed controls) \textit{Let }$\mu
$\textit{\ be an optimal relaxed control minimizing the functional }%
$J$\textit{\ over }$\mathcal{R}$\textit{\ and }$\left(  y_{t}^{\mu},z_{t}%
^{\mu}\right)  $\textit{\ the solution of }$\left(  5\right)  $%
\textit{\ associated with }$\mu$\textit{. }Then, there an unique adapted
process $p^{\mu}$\textit{, which is the solution of the forward stochastic
equation (called adjoint equation),}%
\begin{equation}
\left\{
\begin{array}
[c]{l}%
-dp_{t}^{\mu}=\mathcal{H}_{y}\left(  t,y_{t}^{\mu},z_{t}^{\mu},\mu_{t}%
,p_{t}^{\mu}\right)  dt+\mathcal{H}_{z}\left(  t,y_{t}^{\mu},z_{t}^{\mu}%
,\mu_{t},p_{t}^{\mu}\right)  dW_{t},\\
p_{0}^{\mu}=g_{y}(y_{0}^{\mu}),
\end{array}
\right.
\end{equation}
\textit{such that }%
\begin{equation}
\mathcal{H}\left(  t,y_{t}^{\mu},z_{t}^{\mu},\mu_{t},p_{t}^{\mu}\right)
\geq\mathcal{H}\left(  t,y_{t}^{\mu},z_{t}^{\mu},q_{t},p_{t}^{\mu}\right)
,\ \forall q_{t}\in\mathbb{P}\left(  U\right)  ,\ ae\ ,as.
\end{equation}
where the Hamiltonian $\mathcal{H}$ is defined from $\left[  0,T\right]
\times\mathbb{R}^{n}\times\mathcal{M}_{n\times d}\left(  \mathbb{R}\right)
\times\mathbb{R}^{n}\times\mathbb{P}\left(  U\right)  $\ into $\mathbb{R}$\ by%
\[
\mathcal{H}\left(  t,y,z,p,q\right)  =p%
{\displaystyle\int\nolimits_{U}}
b\left(  t,y,z,a\right)  q_{t}\left(  da\right)  -%
{\displaystyle\int\nolimits_{U}}
h\left(  t,y,z,a\right)  q_{t}\left(  da\right)  .
\]

\end{theorem}

\begin{proof}
Since$\ p_{0}^{\mu}=g_{y}(y_{0}^{\mu})$, then $\left(  23\right)  $ becomes%
\begin{align}
0  &  \leq\mathbb{E}\left[  p_{0}^{\mu}\widetilde{y}_{0}\right]  +\mathbb{E}%
{\displaystyle\int\nolimits_{0}^{T}}
{\displaystyle\int\nolimits_{U}}
\left[  h_{y}\left(  t,y_{t}^{\mu},z_{t}^{\mu},a\right)  \widetilde{y}%
_{t}+h_{z}\left(  t,y_{t}^{\mu},z_{t}^{\mu},a\right)  \widetilde{z}%
_{t}\right]  \mu_{t}\left(  da\right)  dt\\
&  +\mathbb{E}%
{\displaystyle\int\nolimits_{0}^{T}}
\left[
{\displaystyle\int\nolimits_{U}}
h\left(  t,y_{t}^{\mu},z_{t}^{\mu},a\right)  q_{t}\left(  da\right)  -%
{\displaystyle\int\nolimits_{U}}
h\left(  t,y_{t}^{\mu},z_{t}^{\mu},a\right)  \mu_{t}\left(  da\right)
\right]  dt.\nonumber
\end{align}

By applying It\^{o}'s formula to $\left(  p_{t}^{\mu}\widetilde{y}_{t}\right)
$, we have%
\begin{align*}
\mathbb{E}\left[  p_{0}^{\mu}\widetilde{y}_{0}\right]   &  =-\mathbb{E}\left[
%
{\displaystyle\int\nolimits_{0}^{T}}
{\displaystyle\int\nolimits_{U}}
h_{y}\left(  t,y_{t}^{\mu},z_{t}^{\mu},a\right)  \mu_{t}\left(  da\right)
\widetilde{y}_{t}+%
{\displaystyle\int\nolimits_{U}}
h_{z}\left(  t,y_{t}^{\mu},z_{t}^{\mu},a\right)  \mu_{t}\left(  da\right)
\widetilde{z}_{t}\right]  dt\\
&  +\mathbb{E}%
{\displaystyle\int\nolimits_{0}^{T}}
p_{t}^{\mu}\left[
{\displaystyle\int\nolimits_{U}}
b\left(  t,y_{t}^{\mu},z_{t}^{\mu},a\right)  \mu_{t}\left(  da\right)  -%
{\displaystyle\int\nolimits_{U}}
b\left(  t,y_{t}^{\mu},z_{t}^{\mu},a\right)  q_{t}\left(  da\right)  \right]
dt
\end{align*}

Then for every $q\in\mathcal{R}$, $\left(  27\right)  $ becomes%
\begin{equation}
0\leq\mathbb{E}%
{\displaystyle\int\nolimits_{0}^{T}}
\left[  \mathcal{H}\left(  t,y_{t}^{\mu},z_{t}^{\mu},\mu_{t},p_{t}^{\mu
}\right)  -\mathcal{H}\left(  t,y_{t}^{\mu},z_{t}^{\mu},q_{t},p_{t}^{\mu
}\right)  \right]  dt.
\end{equation}

The theorem is proved.
\end{proof}

\subsection{Sufficient optimality conditions for relaxed controls}

In this subsection, we study when necessary optimality conditions $\left(
26\right)  $ becomes sufficient. We recall assumptions $\left(  4\right)  $
and the adjoints equation $\left(  25\right)  $. For any $q\in\mathcal{R}$, we
denote by $\left(  y^{q},z^{q}\right)  $ the solution of equation $\left(
5\right)  $ controlled by $q$.

\begin{theorem}
(Sufficient optimality conditions for relaxed controls).Assume that $g$ and
the function $\left(  y,z\right)  \longmapsto\mathcal{H}\left(
t,y,z,q,p\right)  $ is concave. Then, $\mu$ is an optimal solution of the
relaxed control problem $\left\{  \left(  5\right)  ,\left(  6\right)
,\left(  7\right)  \right\}  $, if it satisfies $\left(  26\right)  .$
\end{theorem}

\begin{proof}
Let $\mu$ be an arbitrary element of $\mathcal{R}$ (candidate to be optimal).
For any $q\in\mathcal{R}$, we have%
\begin{align*}
J\left(  q\right)  -J\left(  \mu\right)   &  =\mathbb{E}\left[  g\left(
y_{0}^{q}\right)  -g\left(  y_{0}^{\mu}\right)  \right] \\
&  +\mathbb{E}%
{\displaystyle\int\nolimits_{0}^{T}}
\left[
{\displaystyle\int\nolimits_{U}}
h\left(  t,y_{t}^{q},z_{t}^{q},a\right)  q_{t}\left(  da\right)  -%
{\displaystyle\int\nolimits_{U}}
h\left(  t,y_{t}^{\mu},z_{t}^{\mu},a\right)  \mu_{t}\left(  da\right)
\right]  dt.
\end{align*}

Since $g$ is convex, then%
\[
g\left(  y_{0}^{q}\right)  -g\left(  y_{0}^{\mu}\right)  \geq g_{y}\left(
y_{0}^{\mu}\right)  \left(  y_{0}^{q}-y_{0}^{\mu}\right)  .
\]

Thus,%
\begin{align*}
J\left(  q\right)  -J\left(  \mu\right)   &  \geq\mathbb{E}\left[
g_{y}\left(  y_{0}^{\mu}\right)  \left(  y_{0}^{q}-y_{0}^{\mu}\right)  \right]
\\
&  +\mathbb{E}%
{\displaystyle\int\nolimits_{0}^{T}}
\left[
{\displaystyle\int\nolimits_{U}}
h\left(  t,y_{t}^{q},z_{t}^{q},a\right)  q_{t}\left(  da\right)  -%
{\displaystyle\int\nolimits_{U}}
h\left(  t,y_{t}^{\mu},z_{t}^{\mu},a\right)  \mu_{t}\left(  da\right)
\right]  dt.
\end{align*}

We remark from $\left(  25\right)  $ that%
\[
p_{0}^{\mu}=g_{y}\left(  y_{0}^{\mu}\right)  .
\]

Then, we have%
\[
J\left(  q\right)  -J\left(  \mu\right)  \geq\mathbb{E}\left[  p_{0}^{\mu
}\left(  y_{0}^{q}-y_{0}^{\mu}\right)  \right]  +\mathbb{E}%
{\displaystyle\int\nolimits_{0}^{T}}
\left[
{\displaystyle\int\nolimits_{U}}
h\left(  t,y_{t}^{q},z_{t}^{q},a\right)  q_{t}\left(  da\right)  -%
{\displaystyle\int\nolimits_{U}}
h\left(  t,y_{t}^{\mu},z_{t}^{\mu},a\right)  \mu_{t}\left(  da\right)
\right]  dt.
\]

By applying It\^{o}'s formula to $p_{t}^{\mu}\left(  y_{t}^{q}-y_{t}^{\mu
}\right)  $, we obtain%
\begin{align*}
\mathbb{E}\left[  p_{0}^{\mu}\left(  y_{0}^{q}-y_{0}^{\mu}\right)  \right]
&  =\mathbb{E}%
{\displaystyle\int\nolimits_{0}^{T}}
\left[  \mathcal{H}_{y}\left(  t,y_{t}^{\mu},z_{t}^{\mu},\mu_{t},p_{t}^{\mu
}\right)  \left(  y_{t}^{q}-y_{t}^{\mu}\right)  +\mathcal{H}_{z}\left(
t,y_{t}^{\mu},z_{t}^{\mu},\mu_{t},p_{t}^{\mu}\right)  \left(  z_{t}^{q}%
-z_{t}^{\mu}\right)  \right]  dt\\
&  +\mathbb{E}%
{\displaystyle\int\nolimits_{0}^{T}}
p_{t}^{\mu}\left[
{\displaystyle\int\nolimits_{U}}
b\left(  t,y_{t}^{\mu},z_{t}^{\mu},a\right)  \mu_{t}\left(  da\right)  -%
{\displaystyle\int\nolimits_{U}}
b\left(  t,y_{t}^{q},z_{t}^{q},a\right)  q_{t}\left(  da\right)  \right]  dt
\end{align*}

Then,%
\begin{align}
&  J\left(  q\right)  -J\left(  \mu\right) \\
&  \geq\mathbb{E}%
{\displaystyle\int\nolimits_{0}^{T}}
\left[  \mathcal{H}\left(  t,y_{t}^{\mu},z_{t}^{\mu},\mu_{t},p_{t}^{\mu
}\right)  -\mathcal{H}\left(  t,y_{t}^{q},z_{t}^{q},q_{t},p_{t}^{\mu}\right)
\right]  dt\nonumber\\
&  +\mathbb{E}%
{\displaystyle\int\nolimits_{0}^{T}}
\mathcal{H}_{y}\left(  t,y_{t}^{\mu},z_{t}^{\mu},\mu_{t},p_{t}^{\mu}\right)
\left(  y_{t}^{q}-y_{t}^{\mu}\right)  dt\nonumber\\
&  +\mathbb{E}%
{\displaystyle\int\nolimits_{0}^{T}}
\mathcal{H}_{z}\left(  t,y_{t}^{\mu},z_{t}^{\mu},\mu_{t},p_{t}^{\mu}\right)
\left(  z_{t}^{q}-z_{t}^{\mu}\right)  dt.\nonumber
\end{align}

Since $\mathcal{H}$ is concave in $\left(  y,z\right)  $ and linear in $\mu$,
then by using the Clarke generalized gradient of $\mathcal{H}$ evaluated at
$\left(  y_{t},z_{t},\mu_{t}\right)  $ and the necessary optimality conditions
$\left(  26\right)  $, it follows by $\left[  33\text{,\ Lemmas\ 2.2 and
2.3}\right]  $ that%
\begin{align*}
\mathcal{H}\left(  t,y_{t}^{\mu},z_{t}^{\mu},\mu_{t},p_{t}^{\mu}\right)
-\mathcal{H}\left(  t,y_{t}^{q},z_{t}^{q},q_{t},p_{t}^{\mu}\right)   &
\geq-\mathcal{H}_{y}\left(  t,y_{t}^{\mu},z_{t}^{\mu},\mu_{t},p_{t}^{\mu
}\right)  \left(  y_{t}^{q}-y_{t}^{\mu}\right) \\
&  -\mathcal{H}_{z}\left(  t,y_{t}^{\mu},z_{t}^{\mu},\mu_{t},p_{t}^{\mu
}\right)  \left(  z_{t}^{q}-z_{t}^{\mu}\right)
\end{align*}

Or equivalently,%
\begin{align*}
0  &  \leq\mathcal{H}\left(  t,y_{t}^{\mu},z_{t}^{\mu},\mu_{t},p_{t}^{\mu
}\right)  -\mathcal{H}\left(  t,y_{t}^{q},z_{t}^{q},q_{t},p_{t}^{\mu}\right)
+\mathcal{H}_{y}\left(  t,y_{t}^{\mu},z_{t}^{\mu},\mu_{t},p_{t}^{\mu}\right)
\left(  y_{t}^{q}-y_{t}^{\mu}\right) \\
&  +\mathcal{H}_{z}\left(  t,y_{t}^{\mu},z_{t}^{\mu},\mu_{t},p_{t}^{\mu
}\right)  \left(  z_{t}^{q}-z_{t}^{\mu}\right)
\end{align*}

Then from $\left(  29\right)  $, we get%
\[
J\left(  q\right)  -J\left(  \mu\right)  \geq0.
\]

The theorem is proved.
\end{proof}

\section{Necessary and sufficient optimality conditions for strict controls}

In this section, we study the strict control problem $\left\{  \left(
1\right)  ,\left(  2\right)  ,\left(  3\right)  \right\}  $ and from the
results of section 3, we derive the optimality conditions for strict controls.

\ 

Throughout this section we suppose moreover that
\begin{align}
&  U\text{ is compact.}\\
&  b\text{ and }h\text{ are\ bounded.}\nonumber
\end{align}

Consider the following subset of $\mathcal{R}$%
\begin{equation}
\delta\left(  \mathcal{U}\right)  =\left\{  q\in\mathcal{R}\text{
\ /\ \ }q=\delta_{v}\ \ ;\ \ v\in\mathcal{U}\right\}  .
\end{equation}

The set $\delta\left(  \mathcal{U}\right)  $ is the collection of all relaxed
controls in the form of Dirac measure charging a strict control.

Denote by $\delta\left(  U\right)  $ the action set of all relaxed controls in
$\delta\left(  \mathcal{U}\right)  $.

If $q\in\delta\left(  \mathcal{U}\right)  $, then $q=\delta_{v}$ with
$v\in\mathcal{U}$. In this case we have for each $t$, $q_{t}\in\delta\left(
U\right)  $ and $q_{t}=\delta_{v_{t}}$.

\ 

We equipped $\mathbb{P}\left(  U\right)  $ with the topology of stable
convergence. Since $U$ is compact, then with this topology $\mathbb{P}\left(
U\right)  $\ is a compact metrizable space. The stable convergence is required
for bounded measurable functions $f\left(  t,a\right)  $\ such that for each
fixed $t\in\left[  0,T\right]  $, $f\left(  t,.\right)  $\ is continuous
(Instead of functions bounded and continuous with respect to the pair $\left(
t,a\right)  $ for the weak topology). The space $\mathbb{P}\left(  U\right)
$\ is equipped with its Borel $\sigma$-field, which is the smallest $\sigma
$-field such that the mapping $q\longmapsto%
{\displaystyle\int}
f\left(  s,a\right)  q\left(  ds,da\right)  $\ are measurable for any bounded
measurable function $f$, continuous with respect to $a.\ $For more details,
see Jacod and Memin $\left[  29\right]  $ and El Karoui et al $\left[
16\right]  $.

This allows us to summarize some of lemmas that we will be used in the sequel.\ 

\ 

\begin{lemma}
(Chattering\thinspace\thinspace Lemma).\thinspace\textit{Let }$q$\textit{\ be
a predictable process with values in the space of probability measures on }%
$U$\textit{. Then there exists a sequence of predictable processes }$\left(
u^{n}\right)  _{n}$\textit{\ with values in }$U$\textit{\ such that }%
\begin{equation}
dtq_{t}^{n}\left(  da\right)  =dt\delta_{u_{t}^{n}}\left(  da\right)
\underset{n\longrightarrow\infty}{\longrightarrow}dtq_{t}\left(  da\right)
\text{ stably},\text{\textit{\ \ }}\mathcal{P}-a.s.
\end{equation}
where $\delta_{u_{t}^{n}}$ is the Dirac measure concentrated at a single point
$u_{t}^{n}$ of $U$.
\end{lemma}

\begin{proof}
See El karoui et al $\left[  14\right]  .$
\end{proof}

\begin{lemma}
Let $q$ be a relaxed control and $\left(  u^{n}\right)  _{n}$ be a sequence of
strict controls such that $\left(  {}\right)  $ holds. Then for any function
$f:\left[  0,T\right]  \times U\rightarrow\mathbb{R}$, measurable in $t$ and
continuous in $a$, we have%
\begin{equation}%
{\displaystyle\int\nolimits_{U}}
f\left(  t,a\right)  \delta_{u_{t}^{n}}\left(  da\right)  \underset
{n\longrightarrow\infty}{\longrightarrow}%
{\displaystyle\int\nolimits_{U}}
f\left(  t,a\right)  q_{t}\left(  da\right)  \ ;\ dt-a.e
\end{equation}

\end{lemma}

\begin{proof}
By $\left(  29\right)  $ and the definition of the stable convergence (See
Jacod-Memin $\left[  21\right]  $), for any bounded function $f\left(
t,a\right)  $ measurable in $t$ and continuous in $a$, we have%
\[%
{\displaystyle\int\nolimits_{0}^{T}}
{\displaystyle\int\nolimits_{U}}
f\left(  t,a\right)  \delta_{u_{t}^{n}}\left(  da\right)  \underset
{n\longrightarrow\infty}{\longrightarrow}%
{\displaystyle\int\nolimits_{0}^{T}}
{\displaystyle\int\nolimits_{U}}
f\left(  t,a\right)  q_{t}\left(  da\right)  .
\]

Put%
\[
g\left(  s,a\right)  =1_{\left[  0,t\right]  }\left(  s\right)  f\left(
s,a\right)  .
\]

It's clear that%
\[%
{\displaystyle\int\nolimits_{0}^{T}}
{\displaystyle\int\nolimits_{U}}
g\left(  s,a\right)  \delta_{u_{s}^{n}}\left(  da\right)  \underset
{n\longrightarrow\infty}{\longrightarrow}%
{\displaystyle\int\nolimits_{0}^{T}}
{\displaystyle\int\nolimits_{U}}
g\left(  s,a\right)  q_{s}\left(  da\right)  .
\]

Then%
\[%
{\displaystyle\int\nolimits_{0}^{t}}
{\displaystyle\int\nolimits_{U}}
f\left(  s,a\right)  \delta_{u_{s}^{n}}\left(  da\right)  \underset
{n\longrightarrow\infty}{\longrightarrow}%
{\displaystyle\int\nolimits_{0}^{t}}
{\displaystyle\int\nolimits_{U}}
f\left(  s,a\right)  q_{s}\left(  da\right)  .
\]

The set $\left\{  \left(  s,t\right)  \ ;\ 0\leq s\leq t\leq T\right\}  $
generate $\mathcal{B}_{\left[  0,T\right]  }$. Then $\forall B\in
\mathcal{B}_{\left[  0,T\right]  }$, we have%
\[%
{\displaystyle\int\nolimits_{B}}
{\displaystyle\int\nolimits_{U}}
f\left(  s,a\right)  \delta_{u_{s}^{n}}\left(  da\right)  \underset
{n\longrightarrow\infty}{\longrightarrow}%
{\displaystyle\int\nolimits_{B}}
{\displaystyle\int\nolimits_{U}}
f\left(  s,a\right)  q_{s}\left(  da\right)  .
\]

This implies that%
\[%
{\displaystyle\int\nolimits_{U}}
f\left(  s,a\right)  \delta_{u_{s}^{n}}\left(  da\right)  \underset
{n\longrightarrow\infty}{\longrightarrow}%
{\displaystyle\int\nolimits_{U}}
f\left(  s,a\right)  q_{s}\left(  da\right)  \ ,\ \ dt-a.e.
\]

The lemma is proved.
\end{proof}

\ \ 

The next lemma gives the stability of the controlled FBSDE with respect to the
control variable.

\begin{lemma}
Let $q\in\mathcal{R}$ be a relaxed control and $\left(  y^{q},z^{q}\right)
$\textit{\ the corresponding trajectory. Then there exists a sequence
}$\left(  u^{n}\right)  _{n}\subset\mathcal{U}$ such that
\begin{align}
\underset{n\rightarrow\infty}{\lim}\mathbb{E}\left[  \underset{t\in\left[
0,T\right]  }{\sup}\left\vert y_{t}^{n}-y_{t}^{q}\right\vert ^{2}\right]   &
=0,\\
\underset{n\rightarrow\infty}{\lim}\mathbb{E}%
{\displaystyle\int\nolimits_{0}^{T}}
\left\vert z_{t}^{n}-z_{t}^{q}\right\vert ^{2}dt  &  =0,
\end{align}%
\begin{equation}
\underset{n\rightarrow\infty}{\lim}J\left(  u^{n}\right)  =J\left(  q\right)
.
\end{equation}
where $\left(  y^{n},z^{n}\right)  $ denotes the solution of equation $\left(
1\right)  $ associated with $u^{n}.$
\end{lemma}

\begin{proof}
i) Proof of $\left(  33\right)  $ and $\left(  34\right)  $.

We have%
\[
\left\{
\begin{array}
[c]{ll}%
d\left(  y_{t}^{n}-y_{t}^{q}\right)  = & -\left[  b\left(  t,y_{t}^{n}%
,z_{t}^{n},u_{t}^{n}\right)  -b\left(  t,y_{t}^{q},z_{t}^{q},u_{t}^{n}\right)
\right]  dt\\
& -\left[  b\left(  t,y_{t}^{q},z_{t}^{q},u_{t}^{n}\right)  -b\left(
t,y_{t}^{q},z_{t}^{q},u_{t}^{n}\right)  \right]  dt\\
& -\left[  b\left(  t,y_{t}^{q},z_{t}^{q},u_{t}^{n}\right)  -%
{\displaystyle\int\nolimits_{U}}
b\left(  t,y_{t}^{q},z_{t}^{q},a\right)  q_{t}\left(  da\right)  \right]  dt\\
& +\left(  z_{t}^{n}-z_{t}^{q}\right)  dW_{t},\\
y_{T}^{n}-y_{T}^{q}= & 0.
\end{array}
\right.
\]

Put%
\begin{align*}
Y_{t}^{n}  &  =y_{t}^{n}-y_{t}^{q},\\
Z_{t}^{n}  &  =z_{t}^{n}-z_{t}^{q},
\end{align*}
and%
\begin{align}
\varphi^{n}\left(  t,Y_{t}^{n},Z_{t}^{n}\right)   &  =b\left(  t,y_{t}%
^{q},z_{t}^{q},u_{t}^{n}\right)  -%
{\displaystyle\int\nolimits_{U}}
b\left(  t,y_{t}^{q},z_{t}^{q},a\right)  q_{t}\left(  da\right) \\
&  +%
{\displaystyle\int\nolimits_{0}^{1}}
b_{y}\left(  t,y_{t}^{q}+\lambda\left(  y_{t}^{n}-y_{t}^{q}\right)  ,z_{t}%
^{q}+\lambda\left(  z_{t}^{n}-z_{t}^{q}\right)  ,u_{t}^{n}\right)  Y_{t}%
^{n}d\lambda\nonumber\\
&  +%
{\displaystyle\int\nolimits_{0}^{1}}
b_{z}\left(  t,y_{t}^{q}+\lambda\left(  y_{t}^{n}-y_{t}^{q}\right)  ,z_{t}%
^{q}+\lambda\left(  z_{t}^{n}-z_{t}^{q}\right)  ,u_{t}^{n}\right)  Z_{t}%
^{n}d\lambda.\nonumber
\end{align}

Then%
\begin{equation}
\left\{
\begin{array}
[c]{l}%
dY_{t}^{n}=\varphi^{n}\left(  t,Y_{t}^{n},Z_{t}^{n}\right)  dt+Z_{t}^{n}%
dW_{t},\\
Y_{T}^{n}=0.
\end{array}
\right.
\end{equation}

The above equation is a linear BSDE with bounded coefficients and with
terminal condition $Y_{T}^{n}=0$. Then by applying a priori estimates (see
Briand et al $\left[  8,\ \text{Proposition 3.2, Page 7}\right]  $), we get%
\[
\mathbb{E}\left[  \underset{t\in\left[  0,T\right]  }{\sup}\left\vert
Y_{t}^{n}\right\vert ^{2}+%
{\displaystyle\int\nolimits_{0}^{T}}
\left\vert Z_{t}^{n}\right\vert ^{2}dt\right]  \leq C\mathbb{E}\left\vert
{\displaystyle\int\nolimits_{0}^{T}}
\left\vert \varphi^{n}\left(  t,0,0\right)  \right\vert dt\right\vert ^{2}.
\]

From $\left(  36\right)  $, we get%
\[
\mathbb{E}\left[  \underset{t\in\left[  0,T\right]  }{\sup}\left\vert
Y_{t}^{n}\right\vert ^{2}+%
{\displaystyle\int\nolimits_{0}^{T}}
\left\vert Z_{t}^{n}\right\vert ^{2}dt\right]  \leq C\mathbb{E}\left\vert
{\displaystyle\int\nolimits_{0}^{T}}
\left\vert b\left(  t,y_{t}^{q},z_{t}^{q},u_{t}^{n}\right)  -%
{\displaystyle\int\nolimits_{U}}
b\left(  t,y_{t}^{q},z_{t}^{q},a\right)  q_{t}\left(  da\right)  \right\vert
dt\right\vert ^{2}.
\]

Since $b$ is continuous and bounded, then from $\left(  32\right)  $ and the
dominated convergence theorem, the term in the right hand side tends to zero
as $n$ tends to infinity. This prove $\left(  33\right)  $ and $\left(
34\right)  $.

\ 

ii) Proof of $\left(  35\right)  .$

Since $g$\ and $h$ are uniformly Lipshitz with respect to $\left(  y,z\right)
$, then by using the Cauchy-Schwartz inequality, we have%
\begin{align*}
&  \left\vert J\left(  q^{n}\right)  -J\left(  q\right)  \right\vert \\
&  \leq C\left(  \mathbb{E}\left\vert y_{0}^{n}-y_{0}^{q}\right\vert
^{2}\right)  ^{1/2}+C\left(  \int_{0}^{T}\mathbb{E}\left\vert y_{t}^{n}%
-y_{t}^{q}\right\vert ^{2}ds\right)  ^{1/2}+C\left(  \mathbb{E}\int_{0}%
^{T}\left\vert z_{t}^{n}-z_{t}^{q}\right\vert ^{2}dt\right)  ^{1/2}\\
&  +C\left(  \mathbb{E}\left\vert \int_{0}^{T}h\left(  t,y_{t}^{q},z_{t}%
^{q},u_{t}^{n}\right)  dt-\int_{0}^{T}\int_{U}h\left(  t,y_{t}^{q},z_{t}%
^{q},a\right)  q_{t}\left(  da\right)  dt\right\vert ^{2}\right)  ^{1/2}.
\end{align*}

From $\left(  33\right)  $ and $\left(  34\right)  $\ the first three terms in
the right hand side converge to zero. Furthermore, since $h$\ is continuous
and bounded, then from $\left(  32\right)  $ and by using the dominated
convergence theorem, the fourth term in the right hand side tends to zero.
This prove $\left(  35\right)  $.
\end{proof}

\begin{lemma}
As a consequence of $\left(  35\right)  $, the strict and the relaxed control
problems have the same value functions. That is
\begin{equation}
\underset{v\in\mathcal{U}}{\inf}J\left(  v\right)  =\underset{q\in\mathcal{R}%
}{\inf}J\left(  q\right)  .
\end{equation}

\end{lemma}

\begin{proof}
Let $u\in\mathcal{U}$ and $\mu\in\mathcal{R}$ be respectively a strict and
relaxed controls such that%
\begin{align}
J\left(  \mu\right)   &  =\underset{q\in\mathcal{R}}{\inf}J\left(  q\right)
.\\
J\left(  u\right)   &  =\underset{v\in\mathcal{U}}{\inf}J\left(  v\right)  ,
\end{align}

By $\left(  39\right)  $, we have
\[
J\left(  \mu\right)  \leq J\left(  q\right)  \text{, }\forall q\in
\mathcal{R}\text{.}%
\]

Since $\delta\left(  \mathcal{U}\right)  \subset\mathcal{R}$, then%
\[
J\left(  \mu\right)  \leq J\left(  q\right)  \text{, }\forall q\in
\delta\left(  \mathcal{U}\right)  \text{.}%
\]

Since $q\in\delta\left(  \mathcal{U}\right)  $, then $q=\delta_{v}$, where
$v\in\mathcal{U}$. It follows%
\[
\left\{
\begin{array}
[c]{c}%
\left(  y^{q},z^{q}\right)  =\left(  y^{v},z^{v}\right)  ,\\
J\left(  q\right)  =J\left(  v\right)  .
\end{array}
\right.
\]

Hence,%
\[
J\left(  \mu\right)  \leq J\left(  v\right)  \text{, }\forall v\in
\mathcal{U}\text{.}%
\]

The control $u$ becomes an element of $\mathcal{U}$, then we get%
\begin{equation}
J\left(  \mu\right)  \leq J\left(  u\right)  \text{.}%
\end{equation}

On the other hand, by $\left(  40\right)  $ we have%
\begin{equation}
J\left(  u\right)  \leq J\left(  v\right)  \text{, }\forall v\in
\mathcal{U}\text{.}%
\end{equation}

The control $\mu$ becomes a relaxed control, then by lemma 12, there exists a
sequence $\left(  v^{n}\right)  _{n}$ of strict controls such that\textit{\ }%
\[
dt\mu_{t}^{n}\left(  da\right)  =dt\delta_{v_{t}^{n}}\left(  da\right)
\underset{n\longrightarrow\infty}{\longrightarrow}dt\mu_{t}\left(  da\right)
\text{ stably},\text{\textit{\ \ }}\mathcal{P}-a.s.
\]

By $\left(  42\right)  $, we get then%
\[
J\left(  u\right)  \leq J\left(  v^{n}\right)  \text{, }\forall n\in
\mathbb{N}\text{,}%
\]

By using $\left(  35\right)  $ and letting $n$ go to infinity in\ the above
inequality, we get%
\begin{equation}
J\left(  u\right)  \leq J\left(  \mu\right)  .
\end{equation}

Finally, by $\left(  41\right)  $ and $\left(  43\right)  $, we have
\[
J\left(  u\right)  =J\left(  \mu\right)  .
\]

The lemma is proved.
\end{proof}

\ 

To establish necessary optimality conditions for strict controls, we need the
following lemma

\ 

\begin{lemma}
The strict control $u$ minimizes $J$ over $\mathcal{U}$ if and only if the
relaxed control $\mu=\delta_{u}$ minimizes $J$ over $\mathcal{R}$.
\end{lemma}

\begin{proof}
Suppose that $u$ minimizes the cost $J$ over $\mathcal{U}$, then
\[
J\left(  u\right)  =\underset{v\in\mathcal{U}}{\inf}J\left(  v\right)
\text{.}%
\]

By $\left(  38\right)  $, we get%
\[
J\left(  u\right)  =\underset{q\in\mathcal{R}}{\inf}J\left(  q\right)
\text{.}%
\]

Since $\mu=\delta_{u}$, then%
\begin{equation}
\left\{
\begin{array}
[c]{c}%
\left(  y^{\mu},z^{\mu}\right)  =\left(  y^{u},z^{u}\right)  ,\\
J\left(  \mu\right)  =J\left(  u\right)  ,
\end{array}
\right.
\end{equation}

This implies that%
\[
J\left(  \mu\right)  =\underset{q\in\mathcal{R}}{\inf J\left(  q\right)  }.
\]

Conversely, if $\mu=\delta_{u}$ minimize $J$ over $\mathcal{R}$, then%
\[
J\left(  \mu\right)  =\underset{q\in\mathcal{R}}{\inf J\left(  q\right)  }.
\]

By $\left(  38\right)  $, we get%
\[
J\left(  \mu\right)  =\underset{v\in\mathcal{U}}{\inf J\left(  v\right)  }.
\]

And By $\left(  44\right)  $, we obtain%
\[
J\left(  u\right)  =\underset{v\in\mathcal{U}}{\inf J\left(  v\right)  }.
\]

The proof is completed.
\end{proof}

\ 

The following lemma, who will be used to establish sufficient optimality
conditions for strict controls, shows that we get the results of the above
lemma if we replace $\mathcal{R}$ by $\mathbb{\delta}\left(  \mathcal{U}%
\right)  .$

\begin{lemma}
The strict control $u$ minimizes $J$ over $\mathcal{U}$ if and only if the
relaxed control $\mu=\delta_{u}$ minimizes $J$ over $\delta\left(
\mathcal{U}\right)  $.
\end{lemma}

\begin{proof}
Let $\mu=\delta_{u}$ be an optimal relaxed control minimizing the cost $J$
over $\delta\left(  \mathcal{U}\right)  $, we have then%
\[
J\left(  \mu\right)  \leq J\left(  q\right)  \text{,\ \ }\forall q\in
\delta\left(  \mathcal{U}\right)  .
\]

Since $q\in\delta\left(  \mathcal{U}\right)  $, then there exists
$v\in\mathcal{U}$ such that $q=\delta_{v}.$

It is easy to see that%
\begin{equation}
\left\{
\begin{array}
[c]{c}%
\left(  y^{\mu},z^{\mu}\right)  =\left(  y^{u},z^{u}\right)  ,\\
\left(  y^{q},z^{q}\right)  =\left(  y^{v},z^{v}\right)  ,\\
J\left(  \mu\right)  =J\left(  u\right)  ,\\
J\left(  q\right)  =J\left(  v\right)  .
\end{array}
\right.
\end{equation}

Then, we get%
\[
J\left(  u\right)  \leq J\left(  v\right)  ,\ \ \forall v\in\mathcal{U}%
\text{.}%
\]

Conversely, let $u$ be a strict control minimizing the cost $J$ over
$\mathcal{U}$. Then%
\[
J\left(  u\right)  \leq J\left(  v\right)  ,\ \ \forall v\in\mathcal{U}%
\text{.}%
\]

Since the controls $u,v$ $\in\mathcal{U}$, then there exist $\mu,q\in
\delta\left(  \mathcal{U}\right)  $ such that
\[
\mu=\delta_{u}\ \ \ ,\ \ \ q=\delta_{v}.
\]

This implies that relations $\left(  45\right)  $ hold. Consequently, we get%
\[
J\left(  \mu\right)  \leq J\left(  q\right)  \text{,\ \ }\forall q\in
\delta\left(  \mathcal{U}\right)  .
\]

The proof is completed.
\end{proof}

\subsection{Necessary optimality conditions for strict controls}

Define the Hamiltonian $H$\ in the strict case from $\left[  0,T\right]
\times\mathbb{R}^{n}\times\mathcal{M}_{n\times d}\left(  \mathbb{R}\right)
\times\mathbb{R}^{n}\times U$\ into $\mathbb{R}$\ by%
\[
H\left(  t,y,z,p,v\right)  =pb\left(  t,y,z,v\right)  -h\left(
t,y,z,v\right)  .
\]

\begin{theorem}
(Necessary optimality conditions for strict controls) \textit{Let }%
$u$\textit{\ be an optimal control minimizing the functional }$J$%
\textit{\ over }$\mathcal{U}$\textit{\ and }$\left(  y_{t}^{u},z_{t}%
^{u}\right)  $\textit{\ the solution of }$\left(  1\right)  $%
\textit{\ associated with }$u$\textit{. }Then there an unique adapted process
$p$\textit{, solution of }%
\begin{equation}
\left\{
\begin{array}
[c]{l}%
dp_{t}^{u}=H_{y}\left(  t,y_{t}^{u},z_{t}^{u},u_{t},p_{t}^{u}\right)
dt+H_{z}\left(  t,y_{t}^{u},z_{t}^{u},u_{t},p_{t}^{u}\right)  dW_{t},\\
p_{0}^{u}=h_{y}(y_{0}^{u}).
\end{array}
\right.
\end{equation}

Such that%
\begin{equation}
H\left(  t,y_{t}^{u},z_{t}^{u},u_{t},p_{t}^{u}\right)  \geq H\left(
t,y_{t}^{u},z_{t}^{u},v_{t},p_{t}^{u}\right)  ;\forall v_{t}\in U;\ ae\ ,\ as.
\end{equation}

\end{theorem}

\begin{proof}
Let $u$ be an optimal solution of the strict control problem $\left\{  \left(
1\right)  ,\left(  2\right)  ,\left(  3\right)  \right\}  $. Then, there exist
$\mu\in\delta\left(  \mathcal{U}\right)  $ such that
\[
\mu=\delta_{u}.
\]

Since $u$ minimizes the cost $J$ over $\mathcal{U}$, then by lemma 15, $\mu$
minimizes $J$ over $\mathcal{R}$. Hence, by the necessary optimality
conditions for relaxed controls (Theorem 10), there exist an unique adapted
process $p^{\mu}$, solution of $\left(  25\right)  $, such that
\begin{equation}
\mathcal{H}\left(  t,y_{t}^{\mu},z_{t}^{\mu},\mu_{t},p_{t}^{\mu}\right)
\geq\mathcal{H}\left(  t,y_{t}^{\mu},z_{t}^{\mu},q_{t},p_{t}^{\mu}\right)
,\ \forall q_{t}\in\mathbb{P}\left(  U\right)  ,\ a.e,\ a.s.
\end{equation}

Since $\mathbb{\delta}\left(  U\right)  \subset\mathbb{P}\left(  U\right)  $,
then we get
\begin{equation}
\mathcal{H}\left(  t,y_{t}^{\mu},z_{t}^{\mu},\mu_{t},p_{t}^{\mu}\right)
\geq\mathcal{H}\left(  t,y_{t}^{\mu},z_{t}^{\mu},q_{t},p_{t}^{\mu}\right)
,\ \forall q_{t}\in\delta\left(  U\right)  ,\ a.e,\ a.s.
\end{equation}

Since $q\in\delta\left(  \mathcal{U}\right)  $, then there exist
$v\in\mathcal{U}$ such that $q=\delta_{v}$.

We note that $v$ is an arbitrary element of $\mathcal{U}$ since $q$ is arbitrary.

Now, since $\mu=\delta_{u}$ and $q=\delta_{v}$, we can easily see that%
\begin{equation}
\left\{
\begin{array}
[c]{c}%
\left(  y^{\mu},z^{\mu}\right)  =\left(  y^{u},z^{u}\right)  ,\\
\left(  y^{q},z^{q}\right)  =\left(  y^{v},z^{v}\right)  ,\\
p^{\mu}=p^{u},\\
\mathcal{H}\left(  t,y_{t}^{\mu},z_{t}^{\mu},\mu_{t},p_{t}^{\mu}\right)
=H\left(  t,y_{t}^{u},z_{t}^{u},u_{t},p_{t}^{u}\right)  ,\\
\mathcal{H}\left(  t,y_{t}^{\mu},z_{t}^{\mu},q_{t},p_{t}^{\mu}\right)
=H\left(  t,y_{t}^{u},z_{t}^{u},v_{t},p_{t}^{u}\right)  ,
\end{array}
\right.
\end{equation}
where $p^{u}$ is the unique solutions of $\left(  46\right)  $.

By using $\left(  49\right)  $ and $\left(  50\right)  $, we can easy deduce
$\left(  47\right)  $. The proof is completed.
\end{proof}

\subsection{Sufficient optimality conditions for strict controls}

We recall assumptions $\left(  4\right)  $ an the adjoint equation $\left(
46\right)  $.

\begin{theorem}
(Sufficient optimality conditions for strict controls) Assume that $g$ is
convex and the function $\left(  y,z\right)  \longmapsto H\left(
t,y,z,q,p\right)  $ is concave. Then, $u$ is an optimal solution of the
control problem $\left\{  \left(  1\right)  ,\left(  2\right)  ,\left(
3\right)  \right\}  $, if it satisfies $\left(  47\right)  .$
\end{theorem}

\begin{proof}
Let $u$ be a strict control (candidate to be optimal) such that necessary
optimality conditions for strict controls hold. That is%
\begin{equation}
H\left(  t,y_{t}^{u},z_{t}^{u},u_{t},p_{t}^{u}\right)  \geq H\left(
t,y_{t}^{u},z_{t}^{u},v_{t},p_{t}^{u}\right)  ,\forall v_{t}\in U,\ a.e,\ a.s.
\end{equation}

The controls $u,v$ are elements of $\mathcal{U}$, then there exist $\mu
,q\in\delta\left(  \mathcal{U}\right)  $ such that%
\[
\mu=\delta_{u}\ \ ,\ \ q=\delta_{v}.
\]

This implies that relations $\left(  50\right)  $ hold, then by $\left(
51\right)  $, we deduce that%
\[
\mathcal{H}\left(  t,y_{t}^{\mu},z_{t}^{\mu},\mu_{t},p_{t}^{\mu}\right)
\geq\mathcal{H}\left(  t,y_{t}^{\mu},z_{t}^{\mu},q_{t},p_{t}^{\mu}\right)
,\ \forall q_{t}\in\mathbb{\delta}\left(  U\right)  ,\ a.e,\ a.s.
\]

Since $H$ is concave in $\left(  y,z\right)  $, it is easy to see that
$\mathcal{H}$ is concave in $\left(  y,z\right)  $, and since $g$ is convex,
then by the same proof that in theorem 11, we show that $\mu$ minimizes the
cost $J$ over $\mathbb{\delta}\left(  \mathcal{U}\right)  $.

Finally by lemma 16, we deduce that $u$ minimizes the cost $J$ over
$\mathcal{U}$. The proof is completed.
\end{proof}

\begin{remark}
The sufficient optimality conditions for strict controls are proved without
assuming neither the convexity of $U$ nor that of $H$ in $v$.
\end{remark}


\begin{thebibliography}{99}                                                                                               %


\bibitem {}S. Bahlali, B. Mezerdi and B. Djehiche, \textit{Approximation and
optimality necessary conditions in relaxed stochastic control problems,}
Journal of Applied Mathematics and Stochastic Analysis, Volume 2006, pp 1-23.

\bibitem {}S. Bahlali and B. Labed,\textit{\ Necessary and sufficient
conditions of optimality for optimal control problem with initial and terminal
costs}, Rand. Operat. and Stoch. Equ, 2006, Vol 14, No3, pp 291-301.

\bibitem {}S. Bahlali, B. Djehiche and B. Mezerdi, \textit{The relaxed maximum
principle in singular control of diffusions}, SIAM J. Control and Optim, 2007,
Vol 46, Issue 2, pp 427-444.

\bibitem {}S. Bahlali, \textit{Necessary and sufficient conditions of
optimality for relaxed and strict control problems}, SIAM J. Control and
Optim, 2008, Vol. 47, No. 4, pp. 2078--2095.

\bibitem {}S.\ Bahlali, \textit{Necessary and sufficient optimality conditions
for relaxed and strict control problems of forward-backward systems},
SIAM\ J.\ on Applied Mathematics, Submitted.

\bibitem {}S. Bahlali, \textit{Necessary and sufficient condition of
optimality for optimal control problem of forward and backward systems},
Theory of Probability and Its Applications ( TVP), In revision.

\bibitem {}A. Bensoussan, \textit{Non linear filtering and stochastic
control}. Proc. Cortona 1981, Lect. notes in Math. 1982, 972, Springer Verlag.

\bibitem {}Ph. Briand, B. Delyon, Y. Hu, E. Pardoux and L. Stoica, $L^{p}$
\textit{Solutions of backward stochastic differential equations}, Sochastic
Process and their Applications, No 108, $2003$, pp 109-129.

\bibitem {}N. Dokuchaev and X. Y. Zhou, \textit{Stochastic controls with
terminal contingent conditions, }Journal Of Mathematical Analysis And
Applications, $1999$, 238, pp 143-165.

\bibitem {}N. El Karoui, N. Huu Nguyen and M. Jeanblanc
Piqu\'{e},\textit{\ Compactification methods in the control of degenerate
diffusions.} Stochastics, Vol. 20, 1987, pp 169-219.

\bibitem {}N. El Karoui and L. Mazliak, \textit{Backward stochastic
differential equations}, $1997$, Addison Wesley, Longman.

\bibitem {}N. El-Karoui, S. Peng, and M. C. Quenez, \textit{Backward
stochastic differential equations in finance}\textbf{,} $1997$, Math. finance 7.

\bibitem {}N. El-Karoui, S. Peng and M.C. Quenez,\textit{\ A dynamic maximum
principle for the optimization of recursive utilities under constraints},
Annals of Applied Probability, 11(2001), pp 664-693.

\bibitem {}W.H. Fleming, \textit{Generalized solutions in optimal stochastic
control}, Differential games and control theory 2, (Kingston conference 1976),
Lect. Notes in Pure and Appl. Math.30, 1978.

\bibitem {}N.F. Framstad, B. Oksendal and A. Sulem, \textit{A sufficient
stochastic maximum principle for optimal control of jump diffusions and
applications to finance}, J. Optim. Theory and applications. 121, $2004$, pp 77-98.

\bibitem {}M. Fuhrman and G. Tessitore, \textit{Existence of optimal
stochastic controls and global solutions of forward-backward stochastic
differential equations},\textit{\ }SIAM J. Control and Optim, $2004$, Vol 43,
N$%
{{}^\circ}%
$ 3, pp 813-830.

\bibitem {}U.G. Haussmann, \textit{General necessary conditions for optimal
control of stochastic systems}, Math. Programming Studies 6, 1976, pp 30-48.

\bibitem {}U.G. Haussmann,\ \textit{A Stochastic maximum principle for optimal
control of diffusions},\ Pitman Research Notes in Math, 1986, Series 151.

\bibitem {}S.Ji and X. Y. Zhou, \textit{A maximum principle for stochastic
optimal control with terminal state constraints, and its applications.
}Commun. Inf. Syst, 2006, 6(4), pp 321-338.

\bibitem {}H.J. Kushner, \textit{Necessary conditions for continuous parameter
stochastic optimization problems}, SIAM J. Control Optim, Vol. 10, 1973, pp 550-565.

\bibitem {}J. Ma and J. Yong, \textit{Solvability of forward-backward SDEs and
the nodal set of Hamilton-Jacobi-Bellman equations}. A Chinese summary appears
in Chinese Ann. Math. Ser. A 16 (1995), no. 4, 532. Chinese Ann. Math. Ser. B
16, 1995, no. 3, pp 279--298.

\bibitem {}J. Ma and J.\ Zhang, \textit{Representation theorems for backward
stochastic differential equations}, Ann. Appl. Probab., $2002$, 12(4), pp
1390-1418, 2002.

\bibitem {}B. Mezerdi and\ S. Bahlali,\ \textit{Approximation in optimal
control of diffusion processes}, Rand. Operat. and Stoch. Equ, 2000, Vol.8, No
4, pp 365-372.

\bibitem {}B. Mezerdi and\ S. Bahlali,\ \textit{Necessary conditions for
optimality in relaxed stochastic control problems},\ Stochastics And Stoch.
Reports, 2002, Vol 73 (3-4), pp 201-218.

\bibitem {}E. Pardoux and S. Peng, \textit{Adapted solutions of backward
stochastic differential equations}, Sys. Control Letters, $1990$, Vol. 14, pp 55-61.

\bibitem {}S. Peng, \textit{A general stochastic maximum principle for optimal
control problems.} SIAM J. Control and Optim.\ $1990$, 28, N${{}^{\circ}}$ 4,
pp 966-979.

\bibitem {}S. Peng, \textit{Backward stochastic differential equations and
application to optimal control}, Appl. Math. Optim. $1993$,\ 27, pp 125-144.

\bibitem {}S. Peng and Z.\ Wu, \textit{Fully coupled forward-backward
stochastic differential equations and applications to optimal control}. SIAM
J. Control Optim., $1999$, 37, no. 3, pp. 825--843.

\bibitem {}J.T.\ Shi and Z. Wu, \textit{The maximum principle for fully
coupled forward-backward stochastic control system}, Acta Automatica Sinica,
Vol 32, No 2, 2006, pp 161-169.

\bibitem {}Z. Wu, \textit{Maximum Principle for Optimal Control Problem of
Fully Coupled Forward-Backward Stochastic Systems}, Systems Sci. Math. Sci,
$1998$, 11, No.3, pp 249-259.

\bibitem {}W. Xu, \textit{Stochastic maximum principle for optimal control
problem of forward and backward system, }J. Austral. Math. Soc. Ser. B 37,
$1995$, pp 172-185.

\bibitem {}J.\ Yong and X.Y. Zhou, \textit{Stochastic controls : Hamilton
systems and HJB\ equations}, vol 43, Springer, New York, 1999.

\bibitem {}X.Y. Zhou, \textit{Sufficient conditions of optimality for
stochastic systems with controllable diffusions}. IEEE Trans. on Automatic
Control, 1996, 41, pp 1176-1179.
\end{thebibliography}
\end{document}